\documentclass{article}
\usepackage{graphicx}
\usepackage[utf8]{inputenc}
\usepackage[T1]{fontenc}
\usepackage{aeguill}
\usepackage{amssymb}
\usepackage[utf8]{inputenc}
\usepackage[T1]{fontenc}
\usepackage{aeguill}
\usepackage{amssymb}
\usepackage{graphicx}
\usepackage[a4paper,pdftex=true]{geometry}

\usepackage{here}
\usepackage{natbib}
\usepackage{bbm}
\usepackage{amsmath}
\usepackage{hyperref}
\usepackage[a4paper,pdftex=true]{geometry}

\newtheorem{theorem}{Theorem}
\newtheorem{lemma}[theorem]{Lemma}

\newtheorem{proposition}[theorem]{Proposition}

\newtheorem{remark}{Remark}

\def\real{\mathbb{R}}
\def\zit{\mathbb{Z}}
\def\nit{\mathbb{N}}
\def\pit{\Bbb{P}}
\def\ee{\Bbb{E}}

\newcommand{\ind}{\mathbf{1}}
\newcommand{\CQFD}{\hfill $\square$}

\newcommand{\Vect}[1]{
  \mathbf{#1}
}
\newcommand{\tr}[1]{
 {#1}^{\!  T}
}
\newcommand{\Mat}[1]{
  \underline{\mathbf{#1}}
}

\newcommand{\ZZ}{\mathbb{Z}}
\newcommand{\RR}{\mathbb{R}}

\allowdisplaybreaks

\title{\bf  Confidence intervals for the Hurst parameter of a fractional Brownian motion 
based on finite sample size
}

\author{Jean-Christophe Breton\footnote{
Laboratoire de Mathématiques, Images et Applications, Université de La Rochelle, France. 
Email: jcbreton@univ-lr.fr}\: 
and
Jean-François Coeurjolly\footnote{
GIPSA-lab and Laboratory Jean Kuntzmann,  Grenoble University, France. 
Email: Jean-Francois.Coeurjolly@upmf-grenoble.fr}
}

\begin{document}

\maketitle

\begin{abstract}
In this paper, we show how concentration inequalities for Gaussian quadratic form can be used 
to propose exact confidence intervals of the Hurst index parametrizing a fractional Brownian motion. 
Both cases where the scaling parameter of the fractional Brownian motion is known or unknown are investigated. 
These intervals are obtained by observing a single discretized sample path of a fractional Brownian motion 
and without any assumption on the parameter~$H$.
\end{abstract}


\tableofcontents

\section{Introduction}

Since the pioneer work of \cite{A-ManVan68}, the fractional Brownian motion (fBm) has become widely popular as well as in a theoretical context as in applications. Fractional Brownian motion can be defined as the only centered Gaussian process, denoted by $(B_H(t))_{t\in \RR}$, with stationary increments and with variance function $v(\cdot)$, 
given by $v(t)=C^2|t|^{2H}$ for all $t\in \RR$. 
The parameter $H \in (0,1)$ (resp. $C>0$) is referred to as the Hurst parameter (resp. the scaling coefficient). 
In particular, when $H=1/2$, it is the standard Brownian motion. 
In general, the fractional Brownian motion is an $H$-self-similar process, that is for all $\delta>0$, $\left( B_H(\delta t) \right)_{t \in \mathbb{R}} \stackrel{d}{=} \delta^H  \left( B_H(t) \right)_{t \in \mathbb{R}}$ (where $\stackrel{d}{=}$ means equal in finite-dimensional distributions) 
with autocovariance function behaving like $O(|k|^{2H-2})$ as $|k|\to +\infty$. 
Thus, the discretized increments of the fractional Brownian motion (called the fractional Gaussian noise) constitute a short-range dependent process, when $H<1/2$, and a long-range dependent process, when $H>1/2$. 
The index $H$ characterizes also the path regularity since the fractal dimension of the fractional Brownian motion is equal to $D=2-H$. General references on self-similar processes and long-memory processes are given in \cite{B-Ber94} or \cite{B-DouOppTaq03}. 

The aim of this paper is to propose confidence intervals for the Hurst parameter based on a single observation of a discretized sample path of the interval $[0,1]$ of a fractional Brownian motion. 
To do so, the most popular strategy consists in using the {\it asymptotic normality} of some estimators of the Hurst parameter, 
see \cite{A-Coe00} for a survey on the estimation of the self-similarity or \cite{A-SheZhuLee07} and \cite{A-Coe08} for more recent discussions in a robust context.
Recently, a new strategy based on concentration inequalities for Gaussian processes obtained by \cite{A-NouVie09} has been proposed by \cite{A-BreNouPec09}. 
In this case, the confidence intervals are {\it non-asymptotic} and they appear to be very interesting when the sample size is moderate.
Our contribution is to improve this direction both from a theoretical and practical point of view. 
In order to present our different contributions, let us first recall the confidence interval proposed by \cite{A-BreNouPec09}.

\begin{proposition} 
\label{prop-IC-BNP}
Assume that one observes a fractional Brownian motion at times $i/n$ for $i=0,\ldots,n+1$ with scaling coefficient $C=1$ and with Hurst parameter satisfying $H\leq H^\star$ for some known $H^\star\in (0,1)$. Fix $\alpha\in (0,1)$, then for all $n$ large enough satisfying  $q_n(\alpha)<(4-4^{H^\star})\sqrt{n}$, where 
$q_n(\alpha):=\frac12\left(b(\alpha)+\sqrt{b(\alpha)^2+852\log\left(\frac2{\alpha}\right)}\right)$ with 
$b(\alpha):= \frac{71}{\sqrt{n}} \log\left(\frac{2}{\alpha}\right)$, we have 
\begin{equation}\label{IC-BNP}
\pit \left( H \in \left[\max\left(0, \widetilde{H}_n^{inf}(q_n(\alpha))\right) , \widetilde{H}_n^{sup}(q_n(\alpha))\right]\right)\geq 1-\alpha,
\end{equation}
where for $t>0$
\begin{eqnarray*}
g_n\left(\widetilde{H}_n^{inf}(t)\right)&:=& \frac12 -\frac{\log(S_n)}{2\log(n)} + \frac{\log\left( 1-\frac{t}{(4-4^{H^\star})\sqrt{n}}\right)}{2\log(n)} \\
g_n\left(\widetilde{H}_n^{sup}(t)\right)&:=&\frac12 -\frac{\log(S_n)}{2\log(n)} + \frac{\log\left( 1+\frac{t}{(4-4^{H^\star})\sqrt{n}}\right)}{2\log(n)}
\end{eqnarray*}
where $g_n$ is the function defined by $g_n(x)=x-\frac{\log(4-4^x)}{2\log(n)}$ 
and $S_n$ is the following statistic 
\begin{equation}
\label{eq-defSn}
S_n:=\frac1{n} \sum_{i=1}^{n} \left(B_H\left(\frac{i+1}n\right)-2B_H\left(\frac{i}n\right)+B_H\left(\frac{i-1}n\right) \right)^2.
\end{equation}
\end{proposition}

Let us give some general comments on this result. 
First, note that this procedure cannot be applied to a fractional Brownian motion whose scaling coefficient $C$ is unknown. 
Secondly, important drawbacks of this procedure rely upon the assumptions made on $H^\star$ and $n$, which exclude the possibility to use this confidence interval when the sample size is small: 
\begin{itemize}
\item Given $\alpha$ and $H^\star$, the following table presents the minimal value of the sample size $n$  in order to ensure that $q_n(\alpha)<(4-4^{H^\star})\sqrt{n}$.
\begin{center}
\begin{tabular}{cccccccccc}
  \hline
& \multicolumn{9}{c}{$H^\star$}\\
 & 0.1 & 0.2 & 0.3 & 0.4 & 0.5 & 0.6 & 0.7 & 0.8 & 0.9 \\
  \hline
$\alpha=1\%$ & 271 & 298 & 335 & 388 & 471 & 611 & 886 & 1592 & 4936 \\
 $\alpha=5\%$ & 189 & 208 & 233 & 270 & 328 & 425 & 617& 1108 & 3437 \\
  $\alpha=10\%$ & 154 & 169 & 190 & 220 & 266 & 346 & 501 & 900 & 2791 \\
   \hline
\end{tabular}
\end{center}
\item The following table exhibits the maximal value of $H^\star$, denoted by $\widetilde{H}^\star$, 
required in order to ensure $q_n(\alpha)<(4-4^{H^\star})\sqrt{n}$ in terms of $\alpha$ and $n$. 
Note that $\widetilde{H}^\star=\log\left( \max\left(1,4-q_n(\alpha)/\sqrt{n}\right) \right)/\log(4)$, 
which means that, given $\alpha$ and $n$, a confidence interval is only available for $H \in (0, \widetilde{H}^\star)$.
\begin{center}
\begin{tabular}{ccccccc}
  \hline
& \multicolumn{5}{c}{$n$} \\
 & 50 & 100 & 200 & 500 & 10000 & 10000 \\
  \hline
$\alpha=1\%$ & 0.00 & 0.00 & 0.00 & 0.53 & 0.93 & 0.93 \\
  $\alpha=5\%$& 0.00 & 0.00 & 0.17 & 0.65 & 0.94 & 0.94 \\
  $\alpha=10\%$ & 0.00 & 0.00 & 0.34 & 0.70 & 0.95 & 0.95 \\
   \hline
\end{tabular}
\end{center}
\end{itemize}
We are now in position to specify our different contributions:
\begin{itemize}
\item We slightly improve the bounds of the concentration inequality obtained by \cite{A-NouVie09}, 
see Section~\ref{sec-ConcIneq} and Proposition ~\ref{prop-CIH2} for more details. 
Note in particular that, in contrast to \cite{A-NouVie09} and \cite{A-BreNouPec09}, we are tracing the constant to optimize numerically our bounds.
\item In the case where the scaling parameter $C$ is known, we propose a new confidence interval without any preliminary assumption on the Hurst parameter $H$ (in contrast to \cite{A-BreNouPec09}) 
and with a very slight condition on the sample size. 
For instance, in comparison to the previous tables, our confidence interval is computable as soon as $n\geq 3$. 
Furthermore, by using ideas similar in \cite{A-Coe01} for the problem of the estimation of the Hurst parameter, 
we also propose a confidence interval when the scaling parameter $C$ is unknown. 
This new confidence interval has the nice property to be independent of $C$ and independent of the discretization step. 
It is remarkable that, in the both cases ($C$ known or unknown), the lengths of the confidence intervals we propose behave 
asymptotically like the ones derived in an asymptotic approach, 
that is they behave like $1/\sqrt{n}\log(n)$ when $C$ is known and $1/\sqrt{n}$ when $C$ is unknown.
\item As suggested by the expression of the statistic in \eqref{eq-defSn},
the procedure described in Proposition~\ref{prop-IC-BNP} is based on the increments of order $2$ of the discretized sample path of the fractional Brownian motion. 
Taking the increments of order $2$ is a special case of filter to work with 
and it is known that discrete filtering has been proposed and used in an estimation context, 
see \cite{A-IstLan97}, \cite{A-KenWoo97}  and \cite{A-Coe01}. 
Recall that the main interest in filtering the fractional Browian motion is that the action of filtering changes the correlation 
so that, for instance, the increments of order $2$ of the fractional Brownian motion constitute a short-range dependent process 
({\it i.e.} its correlation function is absolutely summable). 
Such a  behaviour is required to obtain an efficient concentration inequality. 
In this paper, we propose to construct confidence intervals not only based on the increments of order $2$ 
but on more general filters such as, for instance, increments of larger order or the Daubechies wavelet filters$\ldots$
Finally, let us also underline that a crucial step consists in obtaining an upper-bound of the supremum on the interval $(0,1)$ 
of the $\ell^1-$norm of the correlation function of the discrete filtered series of the fractional Brownian motion. 
When considering the increments of order $2$, \cite{A-BreNouPec09} have obtained the bound $17.75/(4-4^{H^\star})$. 
We have widely improved this point since we compute explicitly this supremum for a large class of filters (including increments of order $2$). 
As an example, for the increments of order~$2$, this gives the explicit value $8/3$.
\item Based on a large simulation study, 
we assess the efficiency of the different procedures that we propose and we compare them with ones based on an asymptotic scheme. 
We discuss and comment these results. 
\end{itemize}

\medskip

The rest of this paper is organized as follows.
In Section \ref{sec-ConcIneq}, we give the concentration inequalities specially designed for our purposes. 
The filtering setting is introduced in Section \ref{sec:applications} where the bounds for the $\ell^1$-norm of the correlation function of the filtered series are also obtained. 
Our confidence intervals for the Hurst parameter are proposed and proved in Section \ref{sec:IC}, both when the scaling parameter is known or unknown. 
Our results are discussed and compared to the literature in Section \ref{sec:simulations}. 
Finally, computations expliciting some bounds for some special filters are given in Appendix \ref{sec:exact}.

\section{Concentration inequalities} 
\label{sec-ConcIneq}

Proposition \ref{prop-IC-BNP} above is based on concentration inequalities proposed by \cite{A-NouVie09} (see Proposition \ref{prop-CI}) 
for smooth enough random variables with respect to Malliavin calculus (see Theorem 4.1-$i)$). 
By applying such inequalities to the random variables $\sqrt{n}V_n$ where $V_n=\frac1n\sum_{i=1}^n H_2(X_i)$, 
$H_2(t)=t^2-1$ is the second Hermite polynomial, 
and $X=\{X_i\}_{1\leq i\leq n}$ is a stationary Gaussian process with variance~$1$ and correlation function $\rho$, 
we obtain concentration inequalities for $H_2-$variations of stationary Gaussian processes. 
In the sequel, for a sequence $(u_i)_{i\in\mathbb{Z}}$, we set $\|u\|_{\ell^1_n}:=\sum_{|i|\leq n} |u_i|$. 
\begin{proposition} 
\label{prop-CIH2} 
Let $\kappa_n =2 \|\rho\|_{\ell^1_n}$. 
Then, for all $t>0$, we have:
\begin{eqnarray}
\label{eq-phir} 
\pit \left( \sqrt{n} V_n \geq t \right) &\leq& \varphi_{r,n}(t;\kappa_n)
:=  e^{-\frac{t\sqrt{n}}{\kappa_n}} \left( 1+ \frac{t}{\sqrt{n}} \right)^{\frac n{\kappa_n}}\\
\label{def-phil}
\pit \left( \sqrt{n} V_n \leq -t \right) &\leq& \varphi_{l,n}(t;\kappa_n)
:=  e^{\frac{t\sqrt{n}}{\kappa_n}} \left( 1- \frac{t}{\sqrt{n}} \right)^{\frac n{\kappa_n}} \ind_{[0,\sqrt{n}]}(t).
\end{eqnarray}
\end{proposition}
Note that Proposition~\ref{prop-CIH2} can be applied to short-memory as well as to long-memory stationary Gaussian processes (as soon as $n$ remains finite).
In order to derive Proposition \ref{prop-CIH2} below, we shall briefly use some notions of Malliavin calculus. 
We just recall the only necessary for our argument and we refer to \cite{A-BreNouPec09} and references therein for any further details. 
We stress that, once Proposition \ref{prop-CIH2} is derived, only basic probability tools will be used. 
Without restriction, we assume the Gaussian random variables $X_i$ have the form $X_i=X(h_i)$ where 
$X(\aleph)=\{X(h): h\in\aleph\}$ is an isonormal Gaussian process over a real separable Hilbert space $\aleph$ 
and $\{h_i:i=0, \dots, n\}$ is a finite subset of $\aleph$ verifying $\ee[X(h_i)X(h_j)]=\rho(i-j)=\langle h_i,h_j\rangle_\aleph$. 
With such a representation, $V_n$ can be seen as a double Wiener-Itô integral with respect to $X$, {\it i.e.} $V_n=I_2\left(\frac 1n\sum_{i=0}^n h_i\otimes h_i\right)$. 
In the sequel, to make easier the presentation, we rewrite Th.~4.1 of \cite{A-NouVie09} only for such random variables, see Proposition~\ref{prop-CI}. 
Actually, in order to optimize our forthcoming results, Proposition~\ref{prop-CI} is a slight improvement of Th.~4.1.
Before, recall that multiple Wiener-Itô integrals $I_q(f)$ are well defined for $f\in\aleph^{\odot q}$, the $q$th symmetric tensor product of $\aleph$, $q\in\nit\setminus\{0\}$;  
the Malliavin derivatives $D$ transforms random variables (in its domain) into random elements with values in~$\aleph$; 
multiple Wiener-Itô integrals are in the domain of $D$ and we have $D_t(I_q(h))=qI_{q-1}(h(\cdot,t))$. 
Recall also that the Hermite polynomials $H_q$ are related to multiple Wiener-Itô integrals by $H_q(I_1(h))=I_q(h^{\otimes q})$ when $\|h\|_\aleph=1$; 
in particular, for $q=2$, we obtain $I_1(h)^2-1=I_2(h^{\otimes 2})$.
\begin{proposition} 
\label{prop-CI}
Let $Z=I_2(f)$ satisfying 
\begin{equation}
\label{eq:condition}
\|DZ\|_{\aleph}^2\leq aZ+b
\end{equation}
for some constants $a\geq 0$ and $b>0$. 
Then, for all $t>0$
\begin{eqnarray*}
\pit(Z\geq t) &\leq& \varphi_{r}(t;a,b) :=e^{-\frac{2t}{a}} \left( 1+ \frac{at}b \right)^{\frac{2b}{a^2}} \\ 
\pit(Z\leq -t) &\leq& \varphi_{l}(t;a,b) :=e^{\frac{2t}{a}} \left( 1- \frac{at}b \right)^{\frac{2b}{a^2}} \mathbf{1}_{[0,b/a]}(t).
\end{eqnarray*}
\end{proposition}
\begin{demo}
The proof is a slight improvement of the bounds in \cite[Theorem 4.1]{A-NouVie09} obtained by a careful reading of the proof
(with the following correspondance with the notation therein: $g_Z(Z)=\frac 12\|DZ\|_{\aleph}^2$, $\alpha=a/2$ and $\beta=b/2$). 
Denoting by $h$ the density of $Z$, the argument of \cite[Theorem 4.1]{A-NouVie09} is based on the following key formula (see (3.16) in \cite{A-NouVie09})
\begin{equation}
\label{eq:key}
\|DZ\|_\aleph^2=\frac{2\int_Z^{+\infty}y h(y)dy}{ h(Z)}.
\end{equation}
For the sake of self-containess, we sketch the main steps of the argument. 
For any $A > 0$, define $m_A : [0,+\infty) \to \real$ by $m_A(\theta)=\ee\left[e^{\theta Z} \ind_{\{Z\leq A\}}\right]$.
We have $m_A'(\theta)=\ee\left[Ze^{\theta Z} \ind_{\{Z\leq A\}}\right]$ and integration by part yields
\begin{eqnarray}
\nonumber
m_A'(\theta)&=&\int_{-\infty}^A xe^{\theta x} h(x)dx\\
\label{eq:maj1}
&\leq& \theta \int_{-\infty}^Ae^{\theta x}\left(\int_x^{+\infty} y h(y)dy\right)dx\\
\label{eq:maj2}
&\leq&\frac\theta 2\ee\left[\|DZ\|_\aleph^2e^{\theta Z}\ind_{\{Z\leq A\}}\right].
\end{eqnarray}
where \eqref{eq:maj1} comes from $\int_A^{+\infty} y h(y)dy\geq 0$ since $\ee[Z]=0$, 
and \eqref{eq:maj2} comes from \eqref{eq:key}.
Because of \eqref{eq:condition}, we obtain for any $\theta\in(0, 2/a)$:
\begin{equation}
\label{eq:ineqd}
m_A'(\theta)\leq \frac{\theta b}{2-\theta a} m_A(\theta).
\end{equation}
Solving \eqref{eq:ineqd}, using $m_A(0)=\pit(Z\leq A)\leq 1$ and applying Fatou's Lemma ($A\to+\infty$) yield
the following bound for the Laplace transform and any $\theta\in(0, 2/a)$:
$$
\ee[e^{\theta Z}]\leq \exp\left(-\frac ba\theta-\frac {2b}{a^2}\ln\left(1-\frac{a\theta}2\right)\right).
$$
The Chebychev inequality together with a standard minimization entail:
$$
\pit(Z\geq t)\leq \exp\left(\min_{\theta\in(0,2/a)}\left\{-\left(t+\frac ba\right)\theta-\frac {2b}{a^2}\ln\left(1-\frac{a\theta}2\right)\right\}\right)
$$
The minimization is achieved in $\widetilde\theta=(2t)/(at+b)$ 
and gives the first bound in Proposition \ref{prop-CI}.
Applying the same argument to $Y=-Z$, satisfying $\|DY\|_\aleph^2\leq -aY+b$, we derive similarly the second bound. 
Note in particular that condition \ref{eq:condition} implies that $Z\geq -b/a$ so that the left tail only makes sense for $t\in (-b/a,0)$. 
\CQFD
\end{demo}

\begin{remark} 
\label{rem-propCI} 
\cite{A-NouVie09} have obtained the bounds 
$$
\phi_l(t;a,b)=\exp\left( -\frac{t^2}{b} \right) \quad \mbox{ and } \quad
\phi_r(t;a,b)=\exp\left( -\frac{t^2}{at+b} \right).
$$
Table~\ref{tab-NVBC} proposes a comparison of these bounds with ours through the comparisons of the values of their reciprocal functions 
since these quantities are of great interest for the considered problem. 
Observe that the most important differences occur when $n$ is moderate. 
The example $a=4/\sqrt{n}$ and $b=4$ corresponds approximately to the choices of parameters that will be used in the next sections.
\end{remark}

\begin{table}[H]
\begin{center}
\begin{tabular}{lrrrrrrrrr}
  \hline
&& \multicolumn{2}{c}{$\alpha=1$\%} & \multicolumn{2}{c}{$\alpha=2.5$\%}&\multicolumn{2}{c}{$\alpha=5$\%} &\multicolumn{2}{c}{$\alpha=10$\%} \\
& & \multicolumn{1}{c}{$\varphi^{-1}_l(\alpha)$}& \multicolumn{1}{c}{$\varphi^{-1}_r(\alpha)$} &\multicolumn{1}{c}{$\varphi^{-1}_l(\alpha)$}& \multicolumn{1}{c}{$\varphi^{-1}_r(\alpha)$} &\multicolumn{1}{c}{$\varphi^{-1}_l(\alpha)$}& \multicolumn{1}{c}{$\varphi^{-1}_r(\alpha)$} &\multicolumn{1}{c}{$\varphi^{-1}_l(\alpha)$}& \multicolumn{1}{c}{$\varphi^{-1}_r(\alpha)$} \\
  \hline
$n=50$&NV & 6.0697 & 9.2102 & 5.4324 & 7.9062 & 4.8955 & 6.8751 & 4.2919 & 5.7878 \\
  &BC &  4.4720 & 7.1547 & 4.1398 & 6.9040 & 3.8372 & 6.0847 & 3.4712 & 5.2008 \\
\hline
 $n=100$ & NV & 6.0697 & 8.1851 & 5.4324 & 7.1048 & 4.8955 & 6.2383 & 4.2919 & 5.3107\\
  &BC &  4.9090 & 7.3551 & 4.4966 & 6.4575 & 4.1314 & 5.7249 & 3.7012 & 4.9267\\
\hline
  $n=500$&NV & 6.0697 & 6.9492 & 5.4324 & 6.1322 & 4.8955 & 5.4606 & 4.2919 & 4.7235\\
  &BC & 5.5334 & 6.6309 & 5.0017 & 5.8810 & 4.5449 & 5.2591 & 4.0218 & 4.5708\\
 \hline
 $n=1000$& NV & 6.0697 & 6.6801 & 5.4324 & 5.9190 & 4.8955 & 5.2891 & 4.2919 & 4.5930 \\
  &BC &   5.6877 & 6.4641 & 5.1259 & 5.7478 & 4.6462 & 5.1513 & 4.1000 & 4.4883\\
 \hline
 $n=10000$ &NV&  6.0697 & 6.2567 & 5.4324 & 5.5819 & 4.8955 & 5.0168 & 4.2919 & 4.3850 \\
  &BC &   5.9475 & 6.1931 & 5.3345 & 5.5312 & 4.8159 & 4.9757 & 4.2308 & 4.3536 \\
   \hline
\end{tabular}
\caption{Computations of the quantities $\varphi^{-1}_l(\alpha)$ and $\varphi^{-1}_r(\alpha)$ for the bounds obtained by \cite{A-NouVie09} (NV) and ours (BC) (see Remark~\ref{rem-propCI} and Proposition~\ref{prop-CI}) for different values of $n$ and $\alpha$ and for the particular case where $a=4/\sqrt{n}$ and $b=4$.}
\label{tab-NVBC}
\end{center}
\end{table}

\begin{remark} 
\label{rem-bij}
Note that $\varphi_{r}(\cdot;a,b)$ (resp. $\varphi_{l}(\cdot;a,b)$) is a bijective function from $(0,+\infty)$ (resp. $(0,b/a)$) to $(0,1)$.
Obviously, the index $l$ in $\varphi_{l}$ (resp. $r$ in $\varphi_{r}$) indicates we consider the left (resp. right) tails.   
\end{remark}

We explain now how Proposition \ref{prop-CIH2} derives from Proposition \ref{prop-CI}:  
standard Malliavin calculus shows that, for $Z=\sqrt n V_n$, $\|DZ\|_\aleph^2=\frac1n\sum_{i,j=1}^n X(i)X(j)\rho(j-i)$, 
see Theorem~2.1 in \cite{A-BreNouPec09}. 
The following lemma ensures that condition \eqref{eq:condition} in Proposition \ref{prop-CI} holds true 
with $a=2\kappa_n/\sqrt{n}$ and $b=2\kappa_n$. 
\begin{lemma}
\label{lemme:ab} 
For $Z=\sqrt n V_n$, we have $\|DZ\|_\aleph^2 \leq \kappa_n\left(\frac 1{\sqrt{n}} Z +1\right)$. 
\end{lemma}
The proof of Lemma \ref{lemme:ab} is a very slight modification of the first part of the proof of Theorem~3.1 in \cite{A-BreNouPec09} to which we refer. 
Finally, Proposition \ref{prop-CI} applies and entails Proposition~\ref{prop-CIH2}.


\section{Applications to quadratic variations of fractional Brownian motion}
\label{sec:applications}

\subsection{Notation}

From now on, $B_H$ stands for a fBm with Hurst parameter $H\in (0,1)$ and with scaling coefficient $C>0$ and 
$\Vect{B}_H$ is the vector of observations at times $i/n$ for $i=0,\ldots,n-1$.  
We consider a filter $\ensuremath{{a^{}}}$ of length $\ell+1$ and order $p$, that is a vector with $\ell+1$ real components $a_i$, $0\leq i\leq \ell$, satisfying 
\begin{equation}
\label{eq:def-filtre}
\sum_{q=0}^\ell q^j a_q=0 \mbox{ for } j=0,\ldots,p-1 \mbox{ and } \sum_{q=0}^\ell q^p a_q \neq 0.
\end{equation}
For instance, we shall consider the following filters: 
Increments 1 ($\ensuremath{{a^{}}}=\{-1,1\}$ with $\ell=1$, $p=1$),
Increments 2 ($\ensuremath{{a^{}}}=\{1,-2,1\}$ with $\ell=2$, $p=2$),
Daublets 4 ($\ensuremath{{a^{}}}=$ $\{-0.09150635$, $-0.15849365$, $0.59150635$, $-0.34150635\}$ with $\ell=3$, $p=2$),
Coiflets 6 ($\ensuremath{{a^{}}}=$ $\{-0.05142973$, $-0.23892973$, $0.60285946$, $-0.27214054$, $-0.05142973$, $0.01107027\}$ with $\ell=5$, $p=2$), 
see {\it e.g.} \cite{A-Dau06} and \cite{B-PerWal00} for more details. 
Let $\Vect{V}^\ensuremath{{a^{}}}$ denote the vector $\Vect{B}_H$ filtered with $\ensuremath{{a^{}}}$ and given for $i=\ell,\ldots,n-1$ by 
$$
V^\ensuremath{{a^{}}}\left(\frac in \right) := \sum_{q=0}^\ell a_q B_H\left(\frac{i-q}n\right).
$$
Let us denote by $\pi_H^{\ensuremath{{a^{}}}}(\cdot)$ and $\rho_H^\ensuremath{{a^{}}}(\cdot)$ the covariance and the correlation functions of the filtered series given by (see \cite{A-Coe01})
\begin{equation}
\label{eq:piH}
\ee[V^\ensuremath{{a^{}}}(k)V^\ensuremath{{a^{}}}(k+j)] = C^2 \times \pi_H^\ensuremath{{a^{}}}(j)
\quad \mbox{ with } \quad
\pi_H^\ensuremath{{a^{}}}(j) = -\frac 12 \sum_{q,r=0}^\ell a_q a_r |q-r+j|^{2H}
\end{equation}
and $\rho_H^\ensuremath{{a^{}}}(\cdot) := \pi_H^\ensuremath{{a^{}}}(\cdot)/\pi_H^\ensuremath{{a^{}}}(0)$ which is independent of $C$. 
Finally, define $S_n^\ensuremath{{a^{}}}$ and $V_n^\ensuremath{{a^{}}}$ as
$$
S_n^\ensuremath{{a^{}}} := \frac{1}{n-\ell}\sum_{i=\ell}^{n-1} V^\ensuremath{{a^{}}}\left( \frac{i}n\right)^2 
$$
and
$$
V_n^\ensuremath{{a^{}}} : = \frac{n^{2H}}{C^2 \pi_H^\ensuremath{{a^{}}}(0)} S_n^\ensuremath{{a^{}}} -1 
= \frac{1}{n-\ell}\sum_{i=\ell}^{n-1} \left(\frac{n^{2H}}{C^2 \pi_H^\ensuremath{{a^{}}}(0)} \times V^\ensuremath{{a^{}}}\left( \frac{i}n\right)^2  -1\right).
$$
Note that $V_n^\ensuremath{{a^{}}}\stackrel{d}{=}\frac1{n-\ell}\sum_{i=\ell}^{n-1}H_2(X^\ensuremath{{a^{}}}_i)$ 
where $H_2(t)=t^2-1$ is the second Hermite polynomial 
and $X^\ensuremath{{a^{}}}$ is a stationary Gaussian process with variance $1$ and with correlation function $\rho_H^\ensuremath{{a^{}}}$. 
Observe that  $V_n^\ensuremath{{a^{}}}$, $n\geq 1$, satisfy a law of large number (LLN) and a central limit theorem (CLT) 
\begin{equation}
\label{eq:LGN_CLT}
V_n^\ensuremath{{a^{}}}\to 0 \mbox{ a.s.}, \quad \sqrt n V_n^\ensuremath{{a^{}}} \Rightarrow {\cal N}(0, \sigma_{H,\ensuremath{{a^{}}}}^2)
\end{equation} 
with explicit variance $\sigma_{H,\ensuremath{{a^{}}}}^2$, see Proposition 1 in \cite{A-Coe01}, used to derive standard confidence interval for $H$. 
In contrast, our argument relies on concentration inequalities: applying Proposition~\ref{prop-CIH2} with these notation, we obtain fo all $s,t\geq 0$:
\begin{equation} 
\label{eq-CIVn}
\pit \left( -s\leq \sqrt{n-\ell} V_n^{\ensuremath{{a^{}}}} \leq t \right) 
\geq 1-\varphi_{r,n-\ell}(t;\kappa_{n,H}^{\ensuremath{{a^{}}}}) -\varphi_{l,n-\ell}(s;\kappa_{n,H}^{\ensuremath{{a^{}}}})
\end{equation}
where $\kappa_{n,H}^a= 2 \sum_{|i|\leq n} |\rho_H^a(i)|$.
As previously explained, the action of filtering a discretized sample path of a fBm changes the correlations into summable correlations for the increments.  
More precisely, it is proved that, for some explicit $k_H$, $\rho_H^\ensuremath{{a^{}}}(i)\sim k_H |i|^{2H-2p}$, see {\it e.g. } \cite{A-Coe01}. 
Thus, $\rho_H^{\ensuremath{{a^{}}}}(\cdot)$ is summable if $p>H+1/2$, 
{\it i.e.} $\rho_H^{\ensuremath{{a^{}}}}(\cdot)$ is summable for all $H\in (0,1)$ for $p\geq 2$ and only for $H\in (0,1/2]$ if $p=1$ 
(in the case $H=1/2$, observe that $\rho_{1/2}^{\ensuremath{{a^{}}}}(k)=0$ for all $|k|\geq \ell$).

One of the aim is to obtain bounds in~\eqref{eq-CIVn} independently of $H$ and easily computable. 
Since $\varphi_{l,n}(t,\cdot)$ and $\varphi_{r,n}(t,\cdot)$ are non-decreasing, the bound \eqref{eq-CIVn} remains true
with $\kappa^a:=2\sup_{H\in (0,\tau)} \|\rho_H^a\|_{\ell^1(\ZZ)}$ replacing $\kappa_{n,H}$. 
Here, and in the sequel, we set $\tau=1/2$ when $p=1$ and $\tau=1$ when $p\geq 2$.  
The following section will prove (among other things) that this quantity is finite.


\subsection{Bounds of $\|\rho_H^{\ensuremath{{a^{}}}}\|_{\ell^1(\mathbb{Z})}$ independent of $H$}

In this section,  we show that $\kappa^{\ensuremath{{a^{}}}}=\sup_{H\in (0,\tau)}  \kappa_H^{\ensuremath{{a^{}}}}$ is finite
for a large class of filters, including the collection of dilated filters $(\ensuremath{{\ensuremath{{a^{}}}^{m}}})_{m\geq 1}$ of a filter $\ensuremath{{a^{}}}$ that will be used in the next section. 
Recall that $\ensuremath{{\ensuremath{{a^{}}}^{m}}}$ is the filter of length $m\ell+1$ with same order $p$ as ${\ensuremath{{a^{}}}}$ 
and defined for $i=0, \ldots, m\ell$ by 
\begin{equation}
\label{def-am}
a_i^m=\left\{ \begin{array}{ll}
a_{i/m} & \mbox{ if } i/m \mbox{ is an integer}\\
0 & \mbox{ otherwise.}
\end{array} \right.
\end{equation}
As a typical example, if $\ensuremath{{a^{}}}:=\ensuremath{{\ensuremath{{a^{}}}^{1}}}=\{1,-2,1\}$, then $\ensuremath{{\ensuremath{{a^{}}}^{2}}}:=\{1,0,-2,0,1\}$.

Since $\pi_H^\ensuremath{{a^{}}}(0)\not=0$, observe that, for a fixed $i \in \ZZ$, 
the functions $H\mapsto \pi_H^\ensuremath{{a^{}}}(i)$ and $H\mapsto \rho_H^\ensuremath{{a^{}}}(i)$ 
are continuous respectively on $[0,1]$ and on $(0,1)$. 
Moreover, since for any filter $\ensuremath{{a^{}}}$,
\begin{equation}
\label{eq:pi0}
\pi_0^\ensuremath{{a^{}}}(0)  = -\frac12 \sum_{q,r=0, q\neq r}^{\ell} a_q a_r = -\frac12 \sum_{q,r=0}^\ell a_q a_r +\frac12 \sum_{q=0}^\ell a_q^2
=\frac12 \sum_{q=0}^\ell a_q^2 >0, 
\end{equation}
the function $H\mapsto \rho_H^\ensuremath{{a^{}}}(i)$ is continuous in 0. 
In particular, this ensures that for $p=1$, $\|\rho_\cdot^{\ensuremath{{a^{}}}}\|_{\ell^1(\ZZ)}$ is continuous on $[0,1/2)$. 
Actually, this may be not continuous in $1/2$ but nevertheless
$\kappa^{\ensuremath{{a^{}}}}=2\sup_{H\in[0,1/2]}\|\rho_H^{\ensuremath{{a^{}}}}\|_{\ell^1(\ZZ)}<+\infty$
for instance $\kappa^{\{-1,1\}}=4$ and $\kappa^{\{-1,1\}^2}=8$. 
We refer to Appendix \ref{sec:exact} for the computation of the exact values 
and to Table \ref{tab-normL1} for the estimation of some other similar constants.

\medskip
For any filter of order $p\geq 2$,  observe that $\pi_1^\ensuremath{{a^{}}}(i)=0$ for all $i$.
Let us consider the following assumption on the filter $\ensuremath{{a^{}}}$, denoted $\mathbf{H}^{\ensuremath{{a^{}}}}$:
\begin{equation}
\label{eq-taua}
\tau^\ensuremath{{a^{}}}:=
\sum_{q,r=0}^{\ell} a_q a_r (q-r)^2\log(|q-r|)  \neq 0,
\end{equation}
with the convention $0\log(0)=0$. 
Tab.~\ref{tab-taua} below shows that Assumption $\mathbf{H^\ensuremath{{a^{}}}}$ is satisfied for a large class of filters. 
Then, from the rule of l'Hospital, 
$$
\lim_{H\to 1^-} \rho_H^{\ensuremath{{a^{}}}}(i) = \frac{ \sum_{q,r=0}^{\ell} a_q a_r (q-r+i)^2\log(|q-r+i|) }{\sum_{q,r=0}^\ell a_q a_r  (q-r)^2\log(|q-r|)}<+\infty.
$$
Therefore, under $\mathbf{H^\ensuremath{{a^{}}}}$, $\rho_H^{\ensuremath{{a^{}}}}(i)$ is a continuous function of $H\in [0,1]$. 
Actually, the same is true for the $\ell^1$-norm of a filter of order $p\geq 2$ as stated in Proposition \ref{prop:CVunif} below. 

\begin{table}[H]
\begin{center}
\begin{tabular}{ll|lllll}
  \hline
& & \multicolumn{5}{|c}{$m$} \\
\multicolumn{2}{c|}{\ensuremath{{a^{}}}} & 1 & 2 & 3 & 4 & 5 \\
\hline
$p=2$&Increments 2 & 5.55 & 22.18 & 49.91 & 88.72 & 138.63 \\
 &Daublets 4  & 0.62 & 2.47 & 5.56 & 9.89 & 15.45 \\
 &Coiflets 6 & 0.61 & 2.42 & 5.45 & 9.69 & 15.15 \\
  $p=3$ &Increments 3 & 13.50 & 53.98 & 121.46 & 215.94 & 337.40 \\
  &Daublets 6 & 0.49 & 1.98 & 4.45 & 7.90 & 12.35 \\
 $p=4$&Increments 4 & 41.43 & 165.70 & 372.84 & 662.82 & 1035.66 \\
 &Daublets 8 & 0.45 & 1.81 & 4.08 & 7.25 & 11.32 \\
 &Symmlets 8 & 0.45 & 1.81 & 4.08 & 7.25 & 11.32 \\
 &Coiflets 12 & 0.45 & 1.79 & 4.03 & 7.16 & 11.19 \\
   \hline
\end{tabular}
\end{center}
\caption{Computations of $\tau^{\ensuremath{{\ensuremath{{a^{}}}^{m}}}}$ for different filters $\ensuremath{{a^{}}}$ 
and its dilatation ${\ensuremath{{a^{}}}^{m}}$ for $m=1,\ldots,5$.} \label{tab-taua}
\end{table}
\begin{proposition}
\label{prop:CVunif}
Let $\ensuremath{{a^{}}}$ be a filter of order $p\geq 2$ satisfying $\mathbf{H}^{\ensuremath{{a^{}}}}$ in \eqref{eq-taua}.
Then $\|\rho_H^{\ensuremath{{a^{}}}}\|_{\ell^1(\ZZ)}$ is a continuous function of $H\in [0,1]$. 
\end{proposition}

\begin{demo}
From \eqref{eq:piH}, we have 
$$
\rho_H^\ensuremath{{a^{}}}(j)
=\frac {|j|^{2H}}{\sum_{q,r=0}^\ell a_qa_r|q-r|^{2H}}\sum_{q,r=0}^\ell a_q a_r \left|1+\frac{q-r}j\right|^{2H}.
$$
For $|j|\geq \ell+1$, we have $q-r+j\geq 0$ for $0\leq q,r\leq\ell$, so that:  
\begin{eqnarray}
\nonumber
\rho_H^\ensuremath{{a^{}}}(j)
&=&\frac {|j|^{2H}}{\sum_{q,r=0}^\ell a_qa_r|q-r|^{2H}}\sum_{q,r=0}^\ell a_q a_r \left(1+\frac{q-r}j\right)^{2H}\\
\nonumber
&=&\frac {|j|^{2H}}{\sum_{q,r=0}^\ell a_qa_r|q-r|^{2H}}\sum_{q,r=0}^\ell a_q a_r \sum_{k=0}^{+\infty} \frac{(2H)(2H-1)\dots(2H-k+1)}{k!}\left(\frac{q-r}j\right)^k\\
\label{eq:outer}
&=&\frac {|j|^{2H}}{\sum_{q,r=0}^\ell a_qa_r|q-r|^{2H}}\sum_{k=2p}^{+\infty} \frac{(2H)(2H-1)\dots(2H-k+1)}{k!j^k}\sum_{q,r=0}^\ell a_q a_r(q-r)^k.
\end{eqnarray}
Observe that in \eqref{eq:outer}, the outer sum starts at $k=2p$. 
This is due to the property \eqref{eq:def-filtre} of the filter $a$ of order $p$ which implies the following remark:
\begin{eqnarray*}
\sum_{q,r=0}^\ell a_q a_r(q-r)^k&=&\sum_{q,r=0}^\ell a_q a_r\sum_{i=0}^k {k\choose i} q^i(-r)^{k-i}\\
&=&\sum_{i=0}^k(-1)^{k-i}{k\choose i}\left(\sum_{q=0}^\ell a_q q^i\sum_{r=0}^\ell a_rr^{k-i}\right)\\
&=&0 \mbox{ if } k\leq 2p-1.
\end{eqnarray*}
As a consequence, for $p\geq 2$, each summand in the outer sum \eqref{eq:outer} contains the factor $2H-2$ in the product $(2H)(2H-1)\dots(2H-k+1)$. 
Observe that under $\mathbf{H}^{\ensuremath{{a^{}}}}$ in \eqref{eq-taua}, 
the rule of l'Hospital ensures that the function $\theta_a(H)=(2-2H)/(\sum_{q\not=r}a_qa_r|q-r|^{2H})$ is bounded at $H=1^-$. 
Since moreover this function is continuous in $H$, we derive, under $\mathbf{H}^{\ensuremath{{a^{}}}}$, 
that $\|\theta_a\|_\infty:=\sup_{H\in[0,1)}|\theta_a(H)|<+\infty$. 

Now, from \eqref{eq:outer}, we have 
\begin{eqnarray}
\nonumber
&&|\rho_H^\ensuremath{{a^{}}}(j)|\\
\nonumber
&=&\left|\theta_a(H)|j|^{2H-2p}\sum_{k=0}^{+\infty} \frac{(2H)(2H-1)(2H-3)\dots(2H-2p-k+1)}{(2p+k)!j^k}\sum_{q,r=0}^\ell a_q a_r(q-r)^{k+2p}\right|\\
\nonumber
&\leq&|\theta_a(H)||j|^{2H-2p}\sum_{k=0}^{+\infty} \frac{(2p+k-1)!}{(2p+k)!j^k}\sum_{q,r=0}^\ell |a_q| |a_r| |q-r|^{k+2p}\\
\nonumber
&\leq&\|\theta_a\|_\infty |j|^{2H-2p}\sum_{q,r=0}^\ell |a_q| |a_r||q-r|^{2p}\sum_{k=0}^{+\infty} \frac1{(k+1)} \left(\frac{|q-r|}{\ell+1}\right)^{k}\\
\label{eq:CVnormale}
&\leq&C(a)|j|^{2H-2p}
\end{eqnarray}
where
\begin{eqnarray*}
C(a)
&=&\|\theta_a\|_\infty\sum_{q,r=0}^\ell |a_q| |a_r||q-r|^{2p} \left(\frac{(\ell+1)\ln(\ell+1)}{\ell}\right) <+\infty. 
\end{eqnarray*}
When $p\geq 2$, the bound \eqref{eq:CVnormale} ensures that the convergence of the series $\sum_{i\in\zit} |\rho_H^a(i)|$ is uniform in $H\in [0,1]$ 
and thus $H\mapsto \|\rho_H^a\|_{\ell^1(\zit)}$ is continuous on $[0,1]$.
\CQFD
\end{demo}

\medskip
Proposition \ref{prop:CVunif} proves the following bound is finite for a filter $\ensuremath{{a^{}}}$ of order $p\geq 2$ satisfying $\mathbf{H^\ensuremath{{a^{}}}}$: 
\begin{equation}
\label{eq:kappa}
\kappa^\ensuremath{{a^{}}}
=2\sup_{H\in(0,1)}\kappa_H^\ensuremath{{a^{}}}
=2\sup_{H\in(0,1)}\|\rho_H^\ensuremath{{a^{}}}\|_{\ell^1(\zit)}
=2\sup_{H\in[0,1]}\|\rho_H^\ensuremath{{a^{}}}\|_{\ell^1(\zit)}<+\infty. 
\end{equation}
As a consequence of this result, this means that the constant $\kappa^{\ensuremath{{a^{}}}}$ 
can be obtained by optimizing the function $H\mapsto\|\rho_H^{\ensuremath{{a^{}}}}\|_{\ell^1(\ZZ)}$ on the interval $[0,1]$.
See Tab. \ref{tab-normL1} below for the computation of such constants for different typical filters. 

For dilated increment-type filters, we manage to compute the exact value of $\|\rho_H^a\|_{\ell^1(\ZZ)}$ (see Appendix~\ref{sec:exact} for more details)
$$\|\rho_H^a\|_{\ell^1(\ZZ)} = 1
+\sum_{k=1}^{\ell-1}\frac{\left|\sum_{j=-\ell}^\ell \alpha_j|j+k|^{2H}\right|}{-\sum_{j=1}^{\ell}\alpha_j j^{2H}}
+(-1)^{p+1}\epsilon(2H-1)\frac{\sum_{k=-\ell+1}^\ell \alpha_k S_{\ell+k-1}^H}{-\sum_{j=1}^\ell\alpha_j j^{2H}},
$$
where $\alpha_j=\sum_{\begin{subarray}{c}q,r=0\\ q-r=j\end{subarray}}^{\ell} a_qa_r$,   $\epsilon(2H-1):=sign(2H-1)$ and where $S_k^H=\sum_{j=0}^{k} j^{2H}$. For the dilated double increments filter $\ensuremath{{a^{}}}=\{1,-2,1\}^m$ for example, this leads to  $\kappa^{\{1,-2,1\}}=2\times 8/3 = 16/3$ and $\kappa^{\{1,-2,1\}^2}= 2\times \left(2 + \frac{25\log(5)-27\log(3)}{8\log(2)} \right) \simeq 7.813554$. 
\begin{table}[H]
\begin{center}
\begin{tabular}{ll|lllll}
 \hline
& & \multicolumn{5}{|c}{$m$} \\
\multicolumn{2}{c|}{\ensuremath{{a^{}}}} & 1 & 2 & 3 & 4 & 5 \\
  \hline
$p=1$&Increments 1 & 2 & 4 & 6 & 8 & 10 \\
$p=2$&Increments 2 & 2.667 & 3.907 & 5.745 & 7.565 & 9.376 \\
&  Daublets 4 & 2.250 & 4.356 & 6.641 & 8.906 & 11.162 \\
&  Coiflets 6 & 2.259 & 4.327 & 6.582 & 8.816 & 11.042 \\
$p=3$&  Increments 3 & 3.200 & 3.783 & 5.396 & 7.406 & 9.200 \\
&  Daublets 6 & 2.429 & 4.516 & 6.688 & 8.833 & 10.966 \\
$p=4$&  Increments 4 & 3.657 & 4.304 & 6.364 & 8.514 & 10.350 \\
&  Daublets 8 & 2.648 & 5.026 & 7.349 & 9.648 & 12.044 \\
&  Coiflets 12 & 2.701 & 5.112 & 7.459 & 9.775 & 12.229 \\
   \hline
\end{tabular}
\end{center} 
\caption{Computation of $\sup_{H\in I} \|\rho_H^\ensuremath{{\ensuremath{{a^{}}}^{m}}}\|_{\ell^1}$ for different filters $\ensuremath{{a^{}}}$ and for $m=1,\ldots,5$. Note that $I=[0,0.5]$ for $p=1$ and $I=[0,1]$ for $p>1$.}\label{tab-normL1}
\end{table}



\section{Confidence intervals of the Hurst parameter}
\label{sec:IC}

For any $\alpha \in (0,1)$, denote by $q_{\bullet,n}^{\ensuremath{{a^{}}}}(\alpha):= \left( \varphi_{\bullet,n}\right)^{-1}(\alpha;\kappa^{\ensuremath{{a^{}}}})$ for $\bullet=l,r$. 
In order to make easier the presentation, define also 
$$
x_{l,n-\ell}^{\ensuremath{{a^{}}}}(\alpha) := 1-\frac{q_{l,n-\ell}^{\ensuremath{{a^{}}}}(\alpha)}{\sqrt{n-\ell}}
\quad \mbox{ and } \quad
x_{r,n-\ell}^{\ensuremath{{a^{}}}}(\alpha) := 1+\frac{q_{r,n-\ell}^{\ensuremath{{a^{}}}}(\alpha)}{\sqrt{n-\ell}}.
$$
Note that Remark~\ref{rem-bij} above ensures that for any $\alpha\in (0,1)$ and for all $n>\ell$, $x_{l,n-\ell}^{\ensuremath{{a^{}}}}(\alpha)>0$.
For further reference, observe that for $\bullet=l,r$ and $n\to+\infty$:
\begin{equation}
\label{eq:q}
q_{\bullet,n-\ell}^{\ensuremath{{a^{}}}}(\alpha)\sim 
q^{\ensuremath{{a^{}}}}(\alpha):=\sqrt{2\kappa^{\ensuremath{{a^{}}}}\log(1/\alpha)}.
\end{equation}
In the sequel, we restrict ourselves, to filters of order $p\geq 2$ which allows us to make no assumption on $H$. Taking a filter of order $p=1$ would have constrained us to assume that $H\leq 1/2$.


\subsection{Scaling parameter $C$ known}

In this section, we assume, without loss of generality, that $C=1$. 
Our confidence interval in Proposition \ref{prop-ICstd} below is expressed in terms of the reciprocal function of 
$g_{n}(x):=2x\log(n)-\log\left( \pi_x^{\ensuremath{{a^{}}}}(0)\right)$, $x\in (0,1)$. 
In order to ensure that $g_{n}$ is indeed invertible, we assume that 
\begin{equation} 
\label{eq-condnICstd}
n\geq \exp\left(\sup_{x\in (0,1)} \frac{ \sum_{q,r=0}^\ell a_q a_r \log(|q-r|) |q-r|^{2x}}{\sum_{q,r=0}^\ell a_q a_r  |q-r|^{2x}}\right).
\end{equation} 
In this case, the function $g_n$ is a strictly increasing bijection from $(0,1)$ to 
$\big( -\log(\pi_0^{\ensuremath{{a^{}}}}(0)), +\infty\big)$.
Moreover recall that a filter of length $\ell+1$ requires a sample size $n\geq \ell+1$. 
Obviously, condition \eqref{eq-condnICstd} only makes sense if the filter $\ensuremath{{a^{}}}$ satisfies:
$$
\sup_{x\in (0,1)} \frac{ \sum_{q,r=0}^\ell a_q a_r \log(|q-r|) |q-r|^{2x}}{\sum_{q,r=0}^\ell a_q a_r  |q-r|^{2x}}<+\infty.
$$
Since $\lim_{x\to 1^-} \sum_{q,r=0}^\ell a_q a_r  |q-r|^{2x} =0^-$
(we stress that this function vanishes with non-positive values of because it is continuous, negative in $x=0$, see \eqref{eq:pi0}, and does not vanish), 
the previous condition is equivalent to the more explicit following one
\begin{equation} 
\label{eq-condaICstd}
\sum_{q,r=0}^\ell a_q a_r \log(|q-r|) (q-r)^2 \geq 0.
\end{equation}
Table \ref{tab:n-eq-condnICstd} exhibits the minimal sample size $n$ required to satisfy~\eqref{eq-condnICstd} 
for different filters $\ensuremath{{\ensuremath{{a^{}}}^{m}}}$ (for $m=1,\ldots,5$) with different order $p=2,3,4$. 
Obviously, condition \eqref{eq-condaICstd} is in force for all these filters. 
\begin{table}
\begin{center}
\begin{tabular}{ll|lllll}
  \hline
&&\multicolumn{5}{|c}{$m$} \\
&$\ensuremath{{a^{}}}$ & 1 & 2 & 3 & 4 & 5 \\
  \hline
$p=2$&Increments 2 & 3 & 4 & 6 & 9 & 11 \\
&Daublets 4& 4 & 6 & 10 & 13 & 15 \\
&Coiflets 6& 6 & 11 & 15 & 21 & 26 \\
$p=3$&Increments 3& 4 & 6 & 10 & 13 & 15 \\
&Daublets 6& 6 & 11 & 15 & 21 & 26 \\
$p=4$&Increments 4& 4 & 9 & 13 & 17 & 21 \\
&Daublets 8& 7 & 15 & 22 & 29 & 36 \\
&Symmlets 8& 7 & 15 & 22 & 29 & 36 \\
&Coiflets 12& 12 & 23 & 34 & 44 & 56 \\
   \hline
\end{tabular}
\end{center}
\caption{Minimal sample size $n$ required to satisfy~\eqref{eq-condnICstd} 
for different dilated filters $\ensuremath{{\ensuremath{{a^{}}}^{m}}}$ of different orders $p$. }
\label{tab:n-eq-condnICstd}
\end{table}

\medskip
\noindent 
We state now our main result when the scaling parameter is known:
\begin{proposition} 
\label{prop-ICstd} 
Let $\alpha \in (0,1)$ be fixed and $\ensuremath{{a^{}}}$ be a filter satisfying $\mathbf{H}^{\ensuremath{{a^{}}}}$ in \eqref{eq-taua}  
\begin{enumerate}
\item For $n\geq \ell+1$, we have:  
\begin{equation}
\label{eq-ICstd-g}
\pit\Big( \log\left(x_{l,n-\ell}^{\ensuremath{{a^{}}}}(\alpha/2)\right) - \log\left( S_n^{\ensuremath{{a^{}}}}\right) \leq g_n(H) \leq \log\left(x_{r,n-\ell}^{\ensuremath{{a^{}}}}(\alpha/2)\right) - \log\left( S_n^{\ensuremath{{a^{}}}}\right) \Big)
\geq 1-\alpha.
\end{equation}
\item Moreover if the filter $\ensuremath{{a^{}}}$ satisfies \eqref{eq-condaICstd} and $n\geq \ell+1$ satisfies \eqref{eq-condnICstd}, we have:  
\begin{equation}
\label{eq-ICstd}
\pit\left( H \in \left[ \widetilde{H}_n^{\inf}(\alpha),\widetilde{H}_n^{\sup}(\alpha)\right] \right) \geq 1-\alpha,
\end{equation}
where
\begin{eqnarray*}
\widetilde{H}_n^{\inf}(\alpha)&:=& \max\Big(0, g_{n}^{-1} \left( \log\left(x_{l,n-\ell}^{\ensuremath{{a^{}}}}(\alpha/2)\right) - \log\left( S_n^{\ensuremath{{a^{}}}}\right) \right) \Big) \\
\widetilde{H}_n^{\sup}(\alpha)&:=& \min \Big(\tau,g_{n}^{-1} \left( \log\left(x_{r,n-\ell}^{\ensuremath{{a^{}}}}(\alpha/2)\right) - \log\left( S_n^{\ensuremath{{a^{}}}}\right) \right) \Big).
\end{eqnarray*}
\item As $n \to +\infty$, the proposed confidence interval in \eqref{eq-ICstd} satisfies almost surely
$$
\left[ \widetilde{H}_n^{\inf}(\alpha),\widetilde{H}_n^{\sup}(\alpha)\right] \rightarrow \{ H\}
$$
and the length $\mu_n$ of the confidence interval satisfies
$$
\mu_n \sim \frac{2 q^{\ensuremath{{a^{}}}}(\alpha/2)}{\sqrt{n}} \frac{1}{g^\prime_n(H)} \sim \frac{q^{\ensuremath{{a^{}}}}(\alpha/2)}{\sqrt{n}\log(n)},
$$
where $q^{\ensuremath{{a^{}}}}$ is defined above in \eqref{eq:q}.
\end{enumerate}
\end{proposition}

\begin{remark}
Proposition~\ref{prop-ICstd} generalizes Proposition~\ref{prop-IC-BNP} derived from \cite{A-BreNouPec09}. 
The scaling parameter is still assumed to be known. However, we do not need to know an upper-bound of $H$ and our condition on $n$ is much sharper than the one required in Proposition~\ref{prop-IC-BNP}. 
As an example, for $\ensuremath{{a^{}}}=(1,-2,1)$, condition~\eqref{eq-condnICstd} is satisfied for all $n \geq 3$, whereas the minimal sample size allowing to derive a confidence interval from Proposition~\ref{prop-IC-BNP} is $1108$ for $\alpha=5\%$ and $H^\star=0.8$. 
\end{remark}

\begin{demo} Consider  the set 
$$
A:= \left\{
 -q_{l,n-\ell}^{\ensuremath{{a^{}}}} (\alpha/2) \leq \sqrt{n-\ell}V_n^{\ensuremath{{a^{}}}} \leq q_{r,n-\ell}^{\ensuremath{{a^{}}}} (\alpha/2)
\right\}.
$$
The bound \eqref{eq-CIVn} entails $\pit(A)\geq 1-\frac\alpha 2 -\frac\alpha 2=1-\alpha$. 
It is now sufficient to notice that
\begin{eqnarray*}
A &=&  \left\{x_{l,n-\ell}^{\ensuremath{{a^{}}}} (\alpha/2)  \leq  1+V_n^{\ensuremath{{a^{}}}} \leq x_{r,n-\ell}^{\ensuremath{{a^{}}}} (\alpha/2) \right\}\\
&=&  \left\{x_{l,n-\ell}^{\ensuremath{{a^{}}}} (\alpha/2)  \leq  \frac{n^{2H}}{\pi_H^{\ensuremath{{a^{}}}}(0)} S_n^{\ensuremath{{a^{}}}} \leq x_{r,n-\ell}^{\ensuremath{{a^{}}}} (\alpha/2) \right\}\\
&=&  \left\{ \log\left(\frac{x_{l,n-\ell}^{\ensuremath{{a^{}}}} (\alpha/2)}{S_n^{\ensuremath{{a^{}}}}}\right)  \leq g_n(H) \leq \log\left(\frac{x_{r,n-\ell}^{\ensuremath{{a^{}}}} (\alpha/2)}{S_n^{\ensuremath{{a^{}}}}} \right) \right\}
\end{eqnarray*}
which proves \eqref{eq-ICstd-g}. 
Next, since under \eqref{eq-condnICstd} and \eqref{eq-condaICstd}, $g_n$ is an increasing bijection, 
\eqref{eq-ICstd} comes immediately from \eqref{eq-ICstd-g}. 
Finally, from \eqref{eq:q}, we have 
\begin{equation*}
\label{eq-logxn}
\log\left(x_{l,n-\ell}^{\ensuremath{{a^{}}}}(\alpha/2)\right) \sim -\frac{q^{\ensuremath{{a^{}}}}(\alpha/2)}{\sqrt{n}} \;\; \mbox{ and } \;\;
\log\left(x_{r,n-\ell}^{\ensuremath{{a^{}}}}(\alpha/2)\right) \sim \frac{q^{\ensuremath{{a^{}}}}(\alpha/2)}{\sqrt{n}}
\end{equation*}
as $n\to+\infty$. 
Moreover, since 
$1+V_n^{\ensuremath{{a^{}}}}= \frac{n^{2H}}{\pi_H^{\ensuremath{{a^{}}}}(0)} S_n^{\ensuremath{{a^{}}}} = S_n^{\ensuremath{{a^{}}}} e^{g_n(H)}$, 
using the LLN in \eqref{eq:LGN_CLT}, we have almost surely
$$
-\log\left( S_n^{\ensuremath{{a^{}}}} \right) = -\log\left( 1+ V_n^{\ensuremath{{a^{}}}} \right) + g_n(H)= g_n(H) -V_n^{\ensuremath{{a^{}}}}(1+o(1)) \sim g_n(H).
$$
It is proved in \cite{A-Coe01} (Proposition~1) that $V_n^a$ converges almost surely towards 0 for any filter and for all $H\in (0,1)$ which implies the almost sure convergence of the confidence interval and the asymptotic behavior of the length $\mu_n$ of the confidence interval.
\CQFD
\end{demo}


\subsection{Scaling parameter $C$ unknown}

The idea to construct confidence intervals when the scaling coefficient $C$ is unknown 
consists in using the collection of the dilated filters \ensuremath{{\ensuremath{{a^{}}}^{m}}} defined in~\eqref{def-am}.  

Let us first introduce some specific notation: let $M\geq 2$ and consider a vector $\Vect{d}=\tr{(d_1,\ldots,d_M)}$ with non zero real components such that $\sum_{i=1}^M d_i=0$
and such that $\tr{\Vect{d}}\Vect{L_M}>0$, where $\Vect{L_M}=\left( \log(m) \right)_{m=1,\ldots,M}$. 
Denote by $I^-$ and $I^+$ the subsets of $\{1,\ldots,M\}$ defined by
$$
I^-=\left\{i\in \{1,\ldots,M\}: d_i <0\right\} 
\quad \mbox{ and } \quad
I^+=\left\{i\in \{1,\ldots,M\}: d_i >0\right\}. 
$$
The following confidence interval is expressed in terms of $\Vect{L_{S_n}} := \left(\log\left(S_n^{\ensuremath{{\ensuremath{{a^{}}}^{m}}}}\right) \right)_{m=1,\ldots,M}$.
\begin{proposition} 
\label{prop-ICgen}
Let $\alpha \in (0,1)$ be fixed and 
denote by $\Vect{L_{X_n^{\inf}}}$ and $\Vect{L_{X_n^{\sup}}}$ the two following vectors with components 
$$
\left(\Vect{L_{X_n^{\inf}}} \right)_m \hskip -5pt= 
\left\{ 
\begin{array}{ll}
\log\left( x^{\ensuremath{{\ensuremath{{a^{}}}^{m}}}}_{l,n-m\ell}\left(\alpha/2M\right) \right) & \!\!\!\!\mbox{if } m\in I^- \\
\log\left( x^{\ensuremath{{\ensuremath{{a^{}}}^{m}}}}_{r,n-m\ell}\left(\alpha/2M\right)\right) & \!\!\!\!\mbox{if } m\in I^+ 
\end{array}\right.
\hskip -5pt,
\left(\Vect{L_{X_n^{\sup}}} \right)_m \hskip -5pt= 
\left\{ 
\begin{array}{ll}
\log\left( x^{\ensuremath{{\ensuremath{{a^{}}}^{m}}}}_{r,n-m\ell}\left(\alpha/2M\right) \right) & \!\!\!\!\mbox{if } m\in I^- \\
\log\left( x^{\ensuremath{{\ensuremath{{a^{}}}^{m}}}}_{l,n-m\ell}\left(\alpha/2M\right)\right) & \!\!\!\!\mbox{if } m\in I^+. 
\end{array}\right.
$$
\begin{enumerate}
\item
Let $n\geq M\ell+1$.  
Then  we have 
\begin{equation} 
\label{eq-ICgen}
\pit\left( H \in \left[\widetilde{H}_{n}^{\inf}(\alpha) , \widetilde{H}_{n}^{\sup}(\alpha) \right] \right)\geq 1-\alpha
\end{equation}
where 
\begin{eqnarray*}
\widetilde{H}_{n}^{\inf}(\alpha) &=& \max\left(0, \frac{1}{2 \tr{\Vect{d}} \Vect{L_M} } \left( \tr{\Vect{d}} \Vect{L_{S_n}} - \tr{\Vect{d}} \Vect{L_{X_n^{\inf}}}  \right) \right) \\
\widetilde{H}_{n}^{\sup}(\alpha) &=& \min\left(1,\frac{1}{2 \tr{\Vect{d}} \Vect{L_M} } \left( \tr{\Vect{d}} \Vect{L_{S_n}} - \tr{\Vect{d}} \Vect{L_{X_n^{\sup}}}  \right)\right).
\end{eqnarray*}
\item
As $n \to +\infty$, the proposed confidence interval in \eqref{eq-ICgen} satisfies almost surely
$$
\left[ \widetilde{H}_n^{\inf}(\alpha),\widetilde{H}_n^{\sup}(\alpha)\right] \rightarrow \{ H\}
$$
and its length $\mu_n$ satisfies
$$
\mu_n:= \frac{ \tr{\Vect{d}} \left( \Vect{L_{X_n^{\inf}}}-\Vect{L_{X_n^{\sup}}} \right)}{2\tr{\Vect{d}} \Vect{L_M} }\sim \frac{1}{\sqrt{n}} \; \frac{\tr{\Vect{d}} \Vect{q_M}(\alpha/2M)}{\tr{\Vect{d}} \Vect{L_M}}
$$
where $\Vect{q_M}(\alpha/2M)$ is the vector of length $M$ with components defined by
$$
\left(\Vect{q_M}(\alpha/2M)\right)_m:= 
\left\{
\begin{array}{rl}
-q^{\ensuremath{{\ensuremath{{a^{}}}^{m}}}}(\alpha/2M) & \mbox{ if } m \in I^-\\
q^{\ensuremath{{\ensuremath{{a^{}}}^{m}}}}(\alpha/2M) & \mbox{ if } m \in I^+
\end{array} \right.
$$
with $q^{\ensuremath{{\ensuremath{{a^{}}}^{m}}}}$ defined in \eqref{eq:q}. 
\end{enumerate}
\end{proposition}

\begin{remark}
Proposition \ref{prop-ICgen} generalizes Proposition~\ref{prop-ICstd} since this new confidence interval does not assume that the scaling parameter, $C$ is known. 
More specifically, note that the definition of the interval does not depend on $C$. 
Note also, that if $\Vect{B_H}$ were not observed on $[0,1)$ but with a dilatation factor, then the confidence interval would remain unchanged.
\end{remark}

\begin{demo}
For $m=1, \ldots, M$, we consider the following event  
\begin{eqnarray*}
A_m&:=&\left\{ x_{l,n-m\ell}^{\ensuremath{{\ensuremath{{a^{}}}^{m}}}}\left(\alpha/2M\right) \leq 1+ V_n^{\ensuremath{{\ensuremath{{a^{}}}^{m}}}} \leq x_{r,n-m\ell}^{\ensuremath{{\ensuremath{{a^{}}}^{m}}}}\left(\alpha/2M \right) \right\}.
\end{eqnarray*}
The bounds \eqref{eq-CIVn} entails that $\pit(A_m) \geq 1 -\frac\alpha{2M} -\frac\alpha{2M} = 1 -\frac\alpha M$. 
First, recall that 
$$
V_n^{\ensuremath{{\ensuremath{{a^{}}}^{m}}}} = \frac{n^{2H}}{C^2 \pi_H^{\ensuremath{{\ensuremath{{a^{}}}^{m}}}}(0)} \; S_n^{\ensuremath{{\ensuremath{{a^{}}}^{m}}}} -1= \gamma \times \frac{1}{m^{2H}} \; S_n^{\ensuremath{{\ensuremath{{a^{}}}^{m}}}} -1\quad \mbox{ with } 
\gamma := \gamma_{C,H,n} = \frac{ n^{2H}}{C^2 \pi_H^\ensuremath{{a^{}}}(0)}.
$$
The crucial point in the definition of the confidence interval relies on the fact that $\gamma$ is independent of $m$. 
Second, note that for $m=1,\ldots,M$: 
\begin{eqnarray*}
&A_m&\\
&=&\left\{ \log\left(x_{l,n-m\ell}^{\ensuremath{{\ensuremath{{a^{}}}^{m}}}}\left(\alpha/2M\right)\right) \leq \log\left(1+ V_n^{\ensuremath{{\ensuremath{{a^{}}}^{m}}}} \right)\leq \log\left(x_{r,n-m\ell}^{\ensuremath{{\ensuremath{{a^{}}}^{m}}}}\left(\alpha/2M \right)\right) \right\} \\
&=& \left\{ \log\left(x_{l,n-m\ell}^{\ensuremath{{\ensuremath{{a^{}}}^{m}}}}(\alpha/2M)\right)  - \log(\gamma) \leq 
\log\left(S_n^{\ensuremath{{\ensuremath{{a^{}}}^{m}}}}\right)- 2H\log(m) \right.\\
&&\hspace*{6cm} \left . \leq \log\left((x_{r,n-m\ell}^{\ensuremath{{\ensuremath{{a^{}}}^{m}}}}(\alpha/2M)\right)  -\log(\gamma)\right\} \\
&=& \left\{ \log\left(S_n^{\ensuremath{{\ensuremath{{a^{}}}^{m}}}}\right)-\log\left(x_{r,n-m\ell}^{\ensuremath{{\ensuremath{{a^{}}}^{m}}}}(\alpha/2M)\right)  + \log(\gamma) \leq 2H\log(m) \right. \\
&&\hspace*{6cm} \left.\leq \log\left(S_n^{\ensuremath{{\ensuremath{{a^{}}}^{m}}}}\right)-\log\left((x_{l,n-m\ell}^{\ensuremath{{\ensuremath{{a^{}}}^{m}}}}(\alpha/2M)\right)+\log(\gamma)
\right\} \\
&=& \left\{ d_m \left(\left(\Vect{L_{S_n}}\right)_m - (\Vect{L_{X_n^{\inf}}})_m + \log(\gamma) \right)\leq 2d_m H (\Vect{L_M})_m \leq 
d_m \left(\left(\Vect{L_{S_n}}\right)_m -  (\Vect{L_{X_n^{\sup}}})_m + \log(\gamma)\right) \right\}.
\end{eqnarray*}
Next, we consider the following event
\begin{eqnarray*}
B&:=& \left\{
\tr{\Vect{d}}\Vect{L_{S_n}} - \tr{\Vect{d}} \Vect{L_{X_n^{\inf}}} +\tr{\Vect{d}}\Vect{1} \log(\gamma)\leq 
2H \tr{\Vect{d}} \Vect{L_M}  \leq 
\tr{\Vect{d}}\Vect{L_{S_n}} - \tr{\Vect{d}} \Vect{L_{X_n^{\sup}}} +\tr{\Vect{d}}\Vect{1} \log(\gamma)
\right\} \\
&=&\left\{
\tr{\Vect{d}}\Vect{L_{S_n}} - \tr{\Vect{d}} \Vect{L_{X_n^{\inf}}} \leq 
2H \tr{\Vect{d}} \Vect{L_M}  \leq
\tr{\Vect{d}}\Vect{L_{S_n}} - \tr{\Vect{d}} \Vect{L_{X_n^{\sup}}} 
\right\} \\
\nonumber
&=& \left\{ H \in \left[ \widetilde{H}_n^{\inf}(\alpha), \widetilde{H}_n^{\sup}(\alpha) \right] \right\}
\end{eqnarray*}
where $\Vect{1}=(1, \dots, 1)^T$. 
Since $A_1 \cap A_2 \cap \ldots \cap A_M \subset B$, setting $A^c=\Omega\setminus A$, we have  
\begin{eqnarray}
\nonumber
\pit(B)& \geq & \pit(A_1 \cap \ldots \cap A_M ) 
= 1-\pit\left((A_1\cap \ldots \cap A_M)^c\right) 
= 1-\pit(A_1^c \cup \ldots \cup A_M^c)  \\
\label{eq:addup}
&\geq & 1 - \sum_{m=1}^M \pit(A_m^c) = \sum_{m=1}^M \pit(A_m) -(M-1)  \\
\nonumber
&\geq& M\left(1-\frac\alpha M \right)-(M-1) = 1-\alpha,
\end{eqnarray}
which ends the proof of \eqref{eq-ICgen}.
Next with the LLN in \eqref{eq:LGN_CLT}, as $n \to +\infty$, the following estimate holds almost surely
\begin{eqnarray*}
\log \left( S_n^{\ensuremath{{\ensuremath{{a^{}}}^{m}}}}\right) &=& 2H \log(m) - \log(\gamma) + \log\left( 1+ V_n^{\ensuremath{{\ensuremath{{a^{}}}^{m}}}}\right) \\
&=& 2H \log(m) - \log(\gamma) + V_n^{\ensuremath{{\ensuremath{{a^{}}}^{m}}}} ( 1+o(1)),
\end{eqnarray*}
and implies that almost surely, when $n\to +\infty$,
\begin{eqnarray*}
\tr{\Vect{d}} \Vect{L_{S_n}} &=& 2H \tr{\Vect{d}} \Vect{L_M} - \tr{\Vect{d}}\mathbf{1} \log(\gamma) + \tr{\Vect{d}} \left( V_n^{\ensuremath{{\ensuremath{{a^{}}}^{m}}}}\right)_{m=1,\ldots,M} (1+o(1)) \\
&=& 2H \tr{\Vect{d}} \Vect{L_M} +  \tr{\Vect{d}} \left( V_n^{\ensuremath{{\ensuremath{{a^{}}}^{m}}}}\right)_{m=1,\ldots,M} (1+o(1)) \\
&\rightarrow& 2H \tr{\Vect{d}} \Vect{L_M} .
\end{eqnarray*}
From~\eqref{eq:q}, one has also the following estimates as $n\to+\infty$:
$$
\left(\Vect{L_{X_n^{\inf}}} \right)_m \sim \frac1{\sqrt{n}} \times
\left\{ 
\begin{array}{rl}
-q^{\ensuremath{{\ensuremath{{a^{}}}^{m}}}}(\alpha/2M) & \!\!\mbox{if } m\in I^- \\
 q^{\ensuremath{{\ensuremath{{a^{}}}^{m}}}}(\alpha/2M)& \!\!\mbox{if } m\in I^+ 
\end{array}\right., \;\;
\left(\Vect{L_{X_n^{\sup}}} \right)_m \sim\frac1{\sqrt{n}} \times
\left\{ 
\begin{array}{rl}
 q^{\ensuremath{{\ensuremath{{a^{}}}^{m}}}}(\alpha/2M)& \!\!\mbox{if } m\in I^- \\
-q^{\ensuremath{{\ensuremath{{a^{}}}^{m}}}}(\alpha/2M)& \!\!\mbox{if } m\in I^+. 
\end{array}\right.
$$
These different results imply the almost sure convergence of the confidence interval towards $\{H\}$. 
For the asymptotic of the length $\mu_n$ of the confidence interval, it is sufficient to note that 
$\left( \Vect{L_{X^{\inf}}} -\Vect{L_{X^{\sup}}} \right) \sim \frac 1{\sqrt n} \Vect{q_M}(\alpha/2M)$. 
\CQFD
\end{demo}


\section{Simulations and discussion}
\label{sec:simulations}

\subsection{Confidence intervals based on the central limit theorem}

\subsubsection{Methodology}

There exists a very wide litterature on the estimation of the Hurst parameter, see {\it e.g.} \cite{A-Coe00} and references therein. 
For all of the available procedures, the confidence interval comes from a limit theorem so that it is of asymptotic very nature. 
In contrast, our confidence intervals in \eqref{eq-ICstd} and \eqref{eq-ICgen} are non-asymptotic since they are based on concentration inequalities. 
In order to compare our procedures, we choose to focus only on one of these procedures which has several similarities with this paper. 
These procedures are based on discrete filtering and are presented in detail in \cite{A-Coe01}. 
For the sake of self-containess, we first summarize them:
\begin{itemize}
\item {\bf Scaling parameter $C$ known.} 
The procedure is based on the fact that almost surely $\frac{n^{2H}}{\pi_H^{\ensuremath{{a^{}}}}(0)}S_n^{\ensuremath{{a^{}}}}\to 1$, $n\to+\infty$. 
With the same function $g_n(x)=2x \log(n) -\log(\pi_x^{\ensuremath{{a^{}}}}(0))$ as the one used to derive the confidence interval in Proposition~\ref{prop-ICstd}, 
this yields the estimator:
$$
\widehat{H}_n^{std}(\ensuremath{{a^{}}}) := g_n^{-1} ( -\log(S_n^\ensuremath{{a^{}}})).
$$
Note that the confidence interval~\eqref{eq-ICstd} is very close to this estimator. 
In particular, the middle of the interval~\eqref{eq-ICstd} behaves asymptotically as $\widehat{H}_n^{std}(\ensuremath{{a^{}}})$.

\item {\bf Scaling parameter $C$ unknown.} 
The idea of \cite{A-Coe00} in this context is to use the following property of quadratic variations of dilated filters $\ee[S_n^{\ensuremath{{\ensuremath{{a^{}}}^{m}}}}]= m^{2H} \gamma$ 
with $\gamma:=\frac{C^2\pi_H^{\ensuremath{{a^{}}}}(0)}{n^{2H}}$ and the almost sure convergence of $S_n^{\ensuremath{{\ensuremath{{a^{}}}^{m}}}}/\ee[S_n^{\ensuremath{{\ensuremath{{a^{}}}^{m}}}}]$ towards 1 for all $m$. 
The idea is then to estimate $H$ via a simple linear regression of $\Vect{L_{S_n}}$ on $2\Vect{L_M}$ for $M$ dilated filters. 
Here, the notation $\Vect{L_{S_n}}$ and $\Vect{L_M}$ are the same as the ones in Proposition~\ref{prop-ICgen}. 
This leads to the estimator
$$
\widehat{H}_n^{gen}(\ensuremath{{a^{}}},M) := \frac{ \tr{\Vect{A}} \Vect{L_{S_n}} }{2\|\Vect{A}\|^2},
$$
where $\Vect{A}=\left(\log(m)-\frac1M\sum_{m=1}^M \log(m)\right)_{m=1, \dots, M}$. 
There is again an analogy between this estimator and our confidence interval in Proposition~\ref{prop-ICgen}. 
Indeed, with $\Vect{d}=\Vect{A}$,  the interval in \eqref{eq-ICgen} rewrites
$$
\left[ \max\left( 0,
\frac{ \tr{\Vect{A}} \left(\Vect{L_{S_n}} - \Vect{L_{X_n^{\inf}}}\right)}{2\|\Vect{A}\|^2}
 \right),
\min\left(1, 
\frac{ \tr{\Vect{A}} \left(\Vect{L_{S_n}} - \Vect{L_{X_n^{\sup}}}\right)}{2\|\Vect{A}\|^2}
\right)
\right],
$$
since $\tr{\Vect{d}}\Vect{L_M}=\tr{\Vect{A}}\Vect{A}=\|A\|^2$.
Again, the middle of this interval behaves asymptotically as $\widehat{H}_n^{gen}(\ensuremath{{a^{}}},M)$. In the particular case $M=2$ the estimator $\widehat{H}_n^{gen}(\ensuremath{{a^{}}},2)$ takes the simple following form
$$
\widehat{H}_n^{gen}(\ensuremath{{a^{}}},2) := \frac{1}{2\log 2} \log\left( \frac{S_n^{\ensuremath{{\ensuremath{{a^{}}}^{2}}}}}{S_n^{\ensuremath{{\ensuremath{{a^{}}}^{1}}}}}  \right)
$$
and the bounds of the interval in \eqref{eq-ICgen} rewrite as
\begin{eqnarray*}
\widetilde{H}_n^{\inf}(\alpha) &:=& \max\left(0, \frac{1}{2\log 2} \left( \log\left( \frac{S_n^{\ensuremath{{\ensuremath{{a^{}}}^{2}}}}}{S_n^{\ensuremath{{\ensuremath{{a^{}}}^{1}}}}}  
\right)
 - \log\left( 
\frac{x_{r,n-2\ell}^{\ensuremath{{\ensuremath{{a^{}}}^{2}}}}(\alpha/4)}{x_{l,n-\ell}^{\ensuremath{{\ensuremath{{a^{}}}^{1}}}}(\alpha/4)}
\right)
\right)\right)
\\
\widetilde{H}_n^{\sup}(\alpha) &:=& \min\left(1, \frac{1}{2\log 2} \left( \log\left( \frac{S_n^{\ensuremath{{\ensuremath{{a^{}}}^{2}}}}}{S_n^{\ensuremath{{\ensuremath{{a^{}}}^{1}}}}}  
\right)
 - \log\left( 
\frac{x_{l,n-2\ell}^{\ensuremath{{\ensuremath{{a^{}}}^{2}}}}(\alpha/4)}{x_{r,n-\ell}^{\ensuremath{{\ensuremath{{a^{}}}^{1}}}}(\alpha/4)}
\right)
\right)\right).
\end{eqnarray*}
\end{itemize}


\subsubsection{Asymptotic confidence intervals}

We refer the reader to \cite{A-Coe01} where the following central limit theorems (CLT) are proved for $\widehat{H}_n^{std}(\ensuremath{{a^{}}})$ and $\widehat{H}_n^{gen}(\ensuremath{{a^{}}},M)$
\begin{equation}
\label{eq-tclstd}
\sqrt{n} \log(n) \frac{\widehat{H}_n^{std}(\ensuremath{{a^{}}}) - H}{\sigma_{std}(\widehat{H}_n^{std})} \stackrel{d}{\longrightarrow} \mathcal{N}(0,1),\quad n\to+\infty
\end{equation}
where $\stackrel{d}{\longrightarrow}$ stands for the convergence in distribution,
$\mathcal{N}(0,1)$ is the normal standard distribution 
and $\sigma^2_{std}(H):=\frac12 \left\|\rho_H^{\ensuremath{{a^{}}}}\right\|_{\ell^2(\ZZ)}$,
and 
\begin{equation} 
\label{eq-tclgen}
\sqrt{n} \frac{\widehat{H}_n^{gen}(\ensuremath{{a^{}}},M) - H}{\sigma_{std}(\widehat{H}_n^{gen},M)} \stackrel{d}{\longrightarrow} \mathcal{N}(0,1), \quad n\to+\infty
\end{equation}
where $\sigma_{gen}^2(H,M):=\frac{\tr{\Vect{A}} \Mat{G} \Vect{A}}{4\|\Vect{A}\|^4}$ 
where $\Mat{G}$ is the $(M\times M)$-matrix defined by $G_{m_1,m_2}= \left\|\rho_H^{\ensuremath{{\ensuremath{{a^{}}}^{m_1}}},\ensuremath{{\ensuremath{{a^{}}}^{m_2}}}}\right\|_{\ell^2(\ZZ)}^2$ 
for $m_1,m_2=1,\ldots,M$, and for all $i \in \ZZ$
$$
\rho_{H}^{{\ensuremath{{\ensuremath{{a^{}}}^{m_1}}},\ensuremath{{\ensuremath{{a^{}}}^{m_2}}}}}(i) =\frac{-\frac12 \sum_{q,r=0}^{\ell} a_a a_r |m_1 q - m_2 r +i|^{2H}  }{\sqrt{\pi_H^{\ensuremath{{\ensuremath{{a^{}}}^{m_1}}}}(0) \pi_H^{\ensuremath{{\ensuremath{{a^{}}}^{m_2}}}}(0)}}.
$$
Note that in the special case where $M=2$, the constant $\sigma^2_{gen}(H,2)$ takes the simple form
$$
\sigma^2_{gen}(H,2)= \frac{1}{2(\log 2)^2} \left( 
\left\| \rho_H^{a^1}\right\|_{\ell^2(\ZZ)}^2
+\left\|\rho_H^{a^2}\right\|_{\ell^2(\ZZ)}^2 
-2 \left\|\rho_H^{a^1,a^2}\right\|_{\ell^2(\ZZ)}^2
\right).
$$
Thanks to the  CLTs, \eqref{eq-tclstd} and \eqref{eq-tclgen} an asymptotic confidence interval to the level $1-\alpha$, $\alpha \in (0,1)$, can be easily constructed
\begin{equation}
\label{eq:ICclt}
IC_\bullet^{clt}(\alpha) =\left[ \max\left(0, \widehat{H}_n^\bullet -\Phi^{-1}(1-\alpha/2) \times \frac{\widehat{\sigma^\bullet}}{v_n^\bullet}\right),
\min\left(1, \widehat{H}_n^\bullet +\Phi^{-1}(1-\alpha/2) \times \frac{\widehat{\sigma^\bullet}}{v_n^\bullet}\right)
\right]
\end{equation}
where $\bullet=std,gen$, $v_n^{std}=\sqrt{n}\log(n)$, $v_n^{gen}=\sqrt{n}$ 
and $\Phi$ is the cumulative distribution function of a standard Gaussian random variable.

\subsection{Comparisons of approaches}

In the following tables, we compare, via Monte-Carlo experiments, the confidence intervals based on concentration inequalities \eqref{eq-ICstd}, \eqref{eq-ICgen} 
and on central limit theorems \eqref{eq:ICclt}. 
The fractional Brownian motions have been generated by using the circulant matrix method ({\it e.g.} \cite{A-KenWoo97}, \cite{A-Coe00}). 
We have realized a very large simulation study. 
The "best" results (in terms of choices of the filters $\ensuremath{{a^{}}}$, of the maximum dilatation factor $M$) are summarized in Table~\ref{tab-simICstd} for the standard fractional Brownian motion ({\it i.e.} $C=1$) 
and in Table~\ref{tab-simICgen} for the general one ({\it i.e.} $C$ unknown).

\begin{figure}[htbp]
\hspace*{-1cm}\begin{tabular}{ll}
\includegraphics[scale=.5]{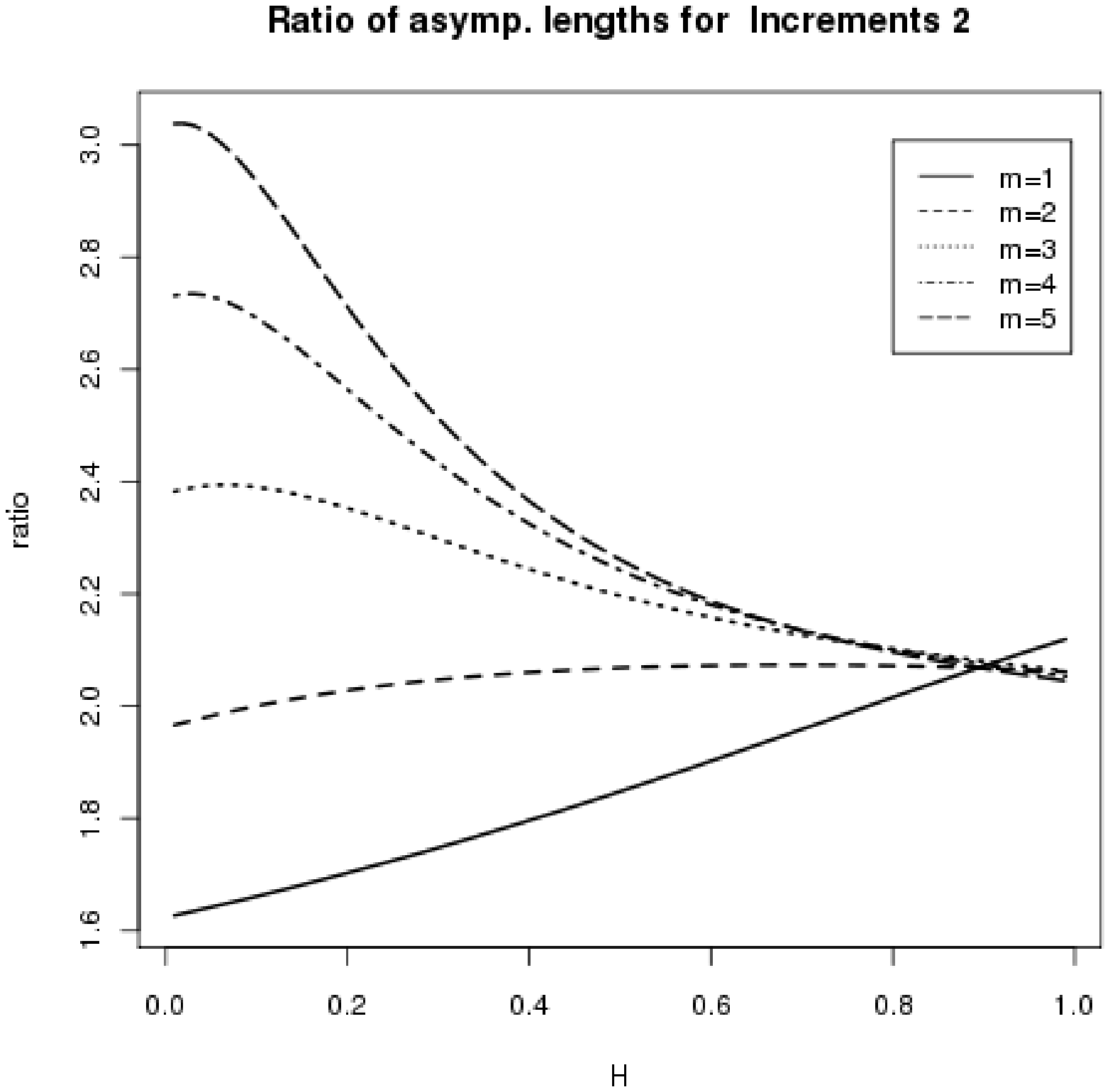} &
\includegraphics[scale=.5]{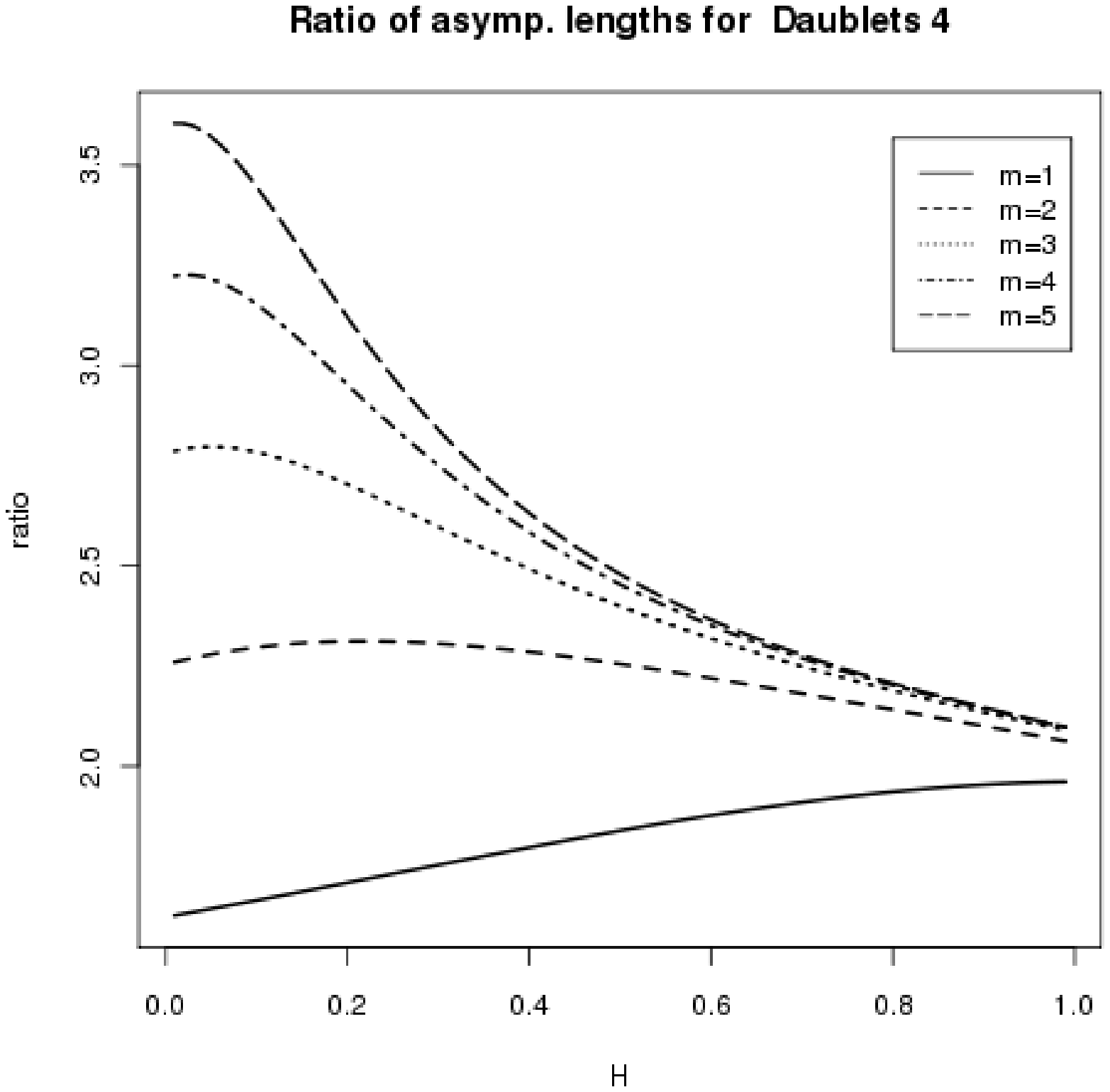} \\
\includegraphics[scale=.5]{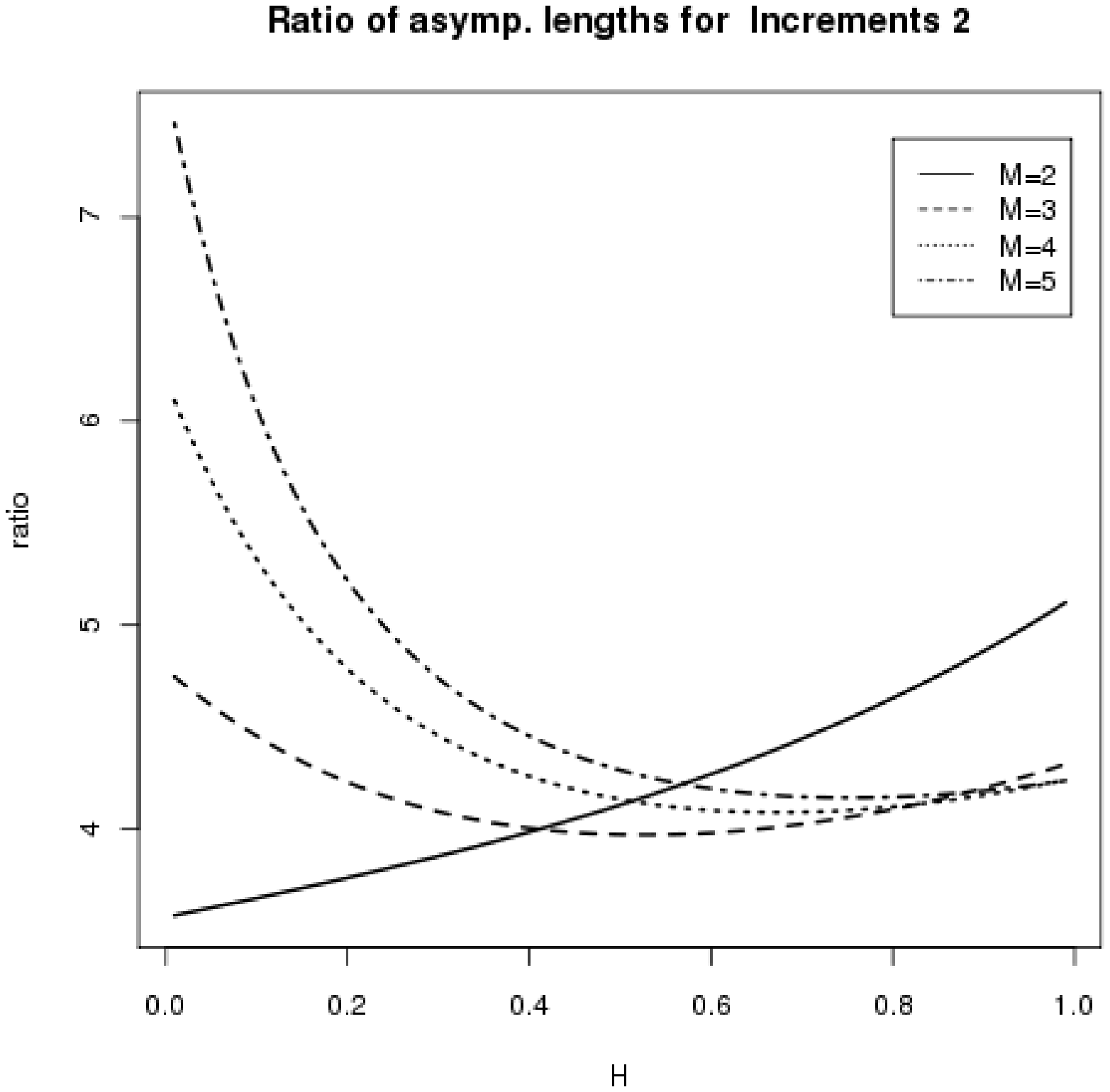} &
\includegraphics[scale=.5]{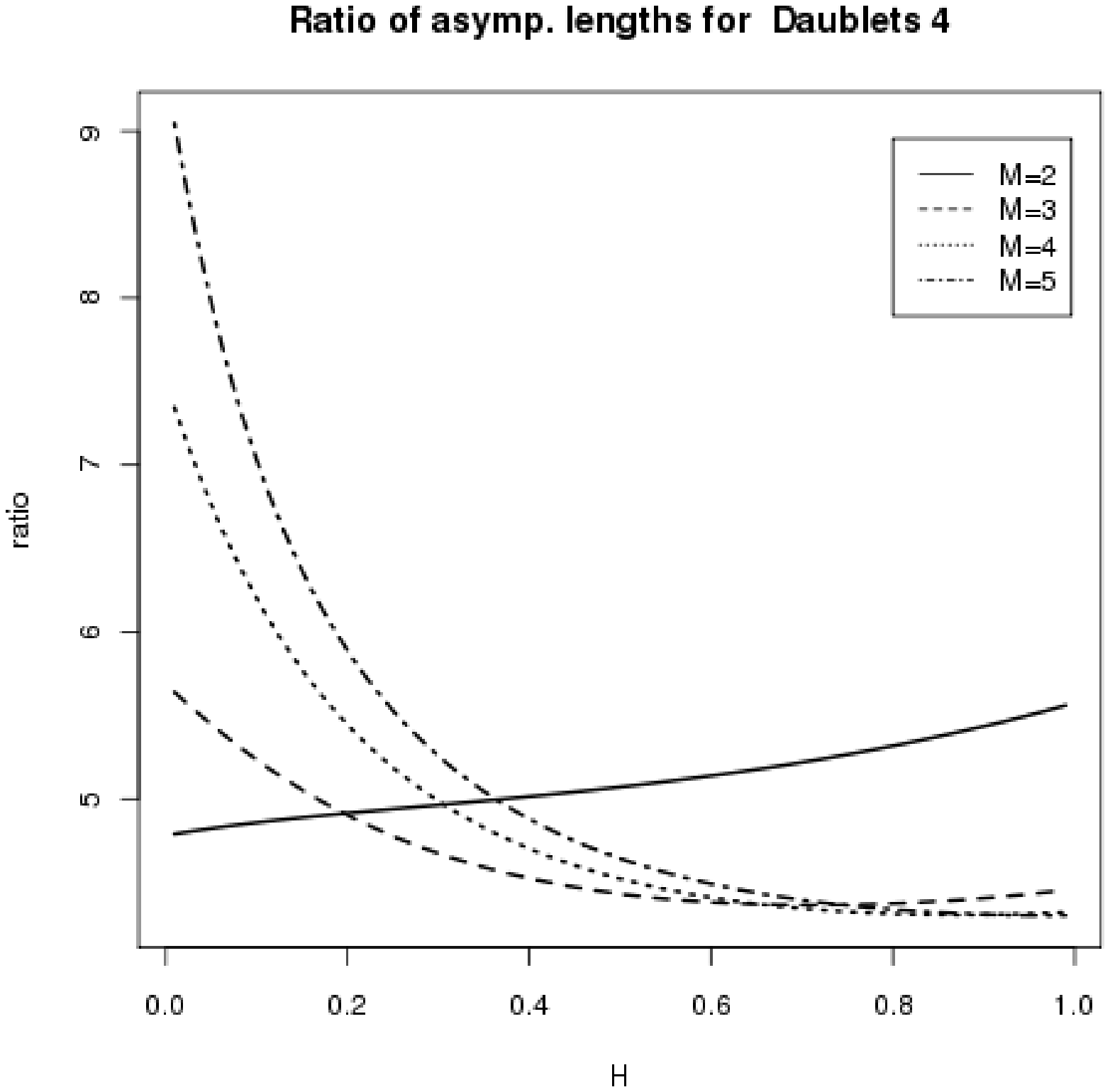} 
\end{tabular}
\caption{Ratio of asymptotic lengths of confidence intervals of procedures derived by concentration inequalities and central limit theorem when the scaling parameter $C$ is known (top) and unknown (bottom). The confidence level equals $1-\alpha=95\%$. For the general procedure, the vector $\Vect{d}$ has been fixed to $\Vect{d}:=\Vect{L_M}-\overline{\Vect{L_M}}$}\label{fig-lengthIC}
\end{figure}

In Figure \ref{fig-lengthIC}, we also compare, in terms of $H$, the asymptotic lengths of the confidence intervals obtained by each approach. 

\begin{table}[H]
\hspace*{-1cm}\begin{tabular}{crrrrrrrrrr}
  \hline
&& \multicolumn{3}{c}{$H=0.2$} & \multicolumn{3}{c}{$H=0.5$} &\multicolumn{3}{c}{$H=0.8$} \\
& & Cover.& Length & $\widehat{H}$ & Cover. & Length & $\widehat{H}$ & Cover. & Length & $\widehat{H}$\\
  \hline
$n=50$&CI[i2] & 100.0 & 0.2191 & 0.1875 & 100.0 & 0.2029 & 0.4832 & 100.0 & 0.1553 & 0.7824 \\
  &CLT[i2] & 95.2 & 0.1330 & 0.2058 & 97.0 & 0.1227 & 0.5013 & 99.6 & 0.1125 & 0.8003 \\
  &CI[d4] & 100.0 & 0.2086 & 0.1886 & 100.0 & 0.1941 & 0.4841 & 100.0 & 0.1482 & 0.7834 \\
  &CLT[d4] & 94.6 & 0.1217 & 0.2050 & 97.2 & 0.1133 & 0.5004 & 99.2 & 0.1076 & 0.7999 \\
\hline
 $n=100$& CI[i2] & 100.0 & 0.1298 & 0.1936 & 100.0 & 0.1212 & 0.4946 & 100.0 & 0.0952 & 0.7931 \\
  & CLT[i2]& 95.0 & 0.0800 & 0.2009 & 97.6 & 0.0737 & 0.5017 & 99.8 & 0.0676 & 0.8003 \\
  &CI[d4] & 100.0 & 0.1224 & 0.1941 & 100.0 & 0.1149 & 0.4949 & 100.0 & 0.0902 & 0.7933 \\
  &CLT[d4] & 95.6 & 0.0732 & 0.2005 & 96.4 & 0.0680 & 0.5012 & 99.6 & 0.0646 & 0.7997 \\
\hline
  $n=500$&CI[i2] & 99.6 & 0.0430 & 0.1994 & 100.0 & 0.0408 & 0.4988 & 99.8 & 0.0336 & 0.7988 \\
  &CLT[i2] & 94.4 & 0.0265 & 0.2004 & 96.4 & 0.0244 & 0.4998 & 98.8 & 0.0224 & 0.7998 \\
  &CI[d4] & 99.6 & 0.0402 & 0.1995 & 100.0 & 0.0383 & 0.4990 & 99.8 & 0.0316 & 0.7989 \\
  &CLT[d4] & 95.4 & 0.0243 & 0.2003 & 96.0 & 0.0225 & 0.4999 & 98.4 & 0.0214 & 0.7998 \\
\hline
  $n=1000$&CI[i2] & 100.0 & 0.0274 & 0.1998 & 100.0 & 0.0262 & 0.4996 & 100.0 & 0.0219 & 0.7997 \\
  &CLT[i2] & 96.6 & 0.0169 & 0.2003 & 97.6 & 0.0155 & 0.5000 & 99.2 & 0.0142 & 0.8001 \\
  &CI[d4] & 100.0 & 0.0256 & 0.1998 & 100.0 & 0.0245 & 0.4996 & 100.0 & 0.0205 & 0.7998 \\
  &CLT[d4] & 96.4 & 0.0154 & 0.2002 & 97.2 & 0.0143 & 0.5000 & 98.8 & 0.0136 & 0.8001 \\
\hline
  $n=10000$&CI[i2] & 99.8 & 0.0066 & 0.2000 & 100.0 & 0.0063 & 0.4999 & 100.0 & 0.0055 & 0.8000 \\
  &CLT[i2] & 94.2 & 0.0040 & 0.2000 & 96.2 & 0.0037 & 0.5000 & 98.4 & 0.0034 & 0.8000 \\
  &CI[d4] & 99.8 & 0.0061 & 0.2000 & 99.8 & 0.0059 & 0.5000 & 100.0 & 0.0051 & 0.8000 \\
  &CLT[d4]& 94.4 & 0.0037 & 0.2000 & 95.0 & 0.0034 & 0.5000 & 98.2 & 0.0032 & 0.8000 \\
   \hline
\end{tabular}
\caption{Monte-carlo experiments based on $500$ replications of a fractional Brownian motion with Hurst parameter $H=0.2,0.5,0.8$ and scaling coefficient $C=1$ (assumed to be known) and for different values of the sample size $n$. 
The filters i2 and d4 denote respectively the filter of Increments of order 2 and the Daublets 4.} \label{tab-simICstd}
\end{table}

\begin{table}[H]
\begin{center}
\begin{tabular}{crrrrrrrrrr}
  \hline
&& \multicolumn{3}{c}{$H=0.2$} & \multicolumn{3}{c}{$H=0.5$} &\multicolumn{3}{c}{$H=0.8$} \\
& & Cover.& Length & $\widehat{H}$ & Cover. & Length & $\widehat{H}$ & Cover. & Length & $\widehat{H}$\\
  \hline
$n=50$&CLT[i2,2] & 95.4 & 0.5970 & 0.3225 & 92.2 & 0.6776 & 0.5064 & 97.2 & 0.5422 & 0.7062 \\
  &CI[i2,2] & 100.0 & 1.0000 & 0.5000 & 100.0 & 1.0000 & 0.5000 & 100.0 & 1.0000 & 0.5000 \\
  &CLT[i2,5] & 89.4 & 0.3706 & 0.2121 & 88.2 & 0.5083 & 0.4838 & 94.2 & 0.4595 & 0.7265 \\
  &CI[i2,5] & 100.0 & 1.0000 & 0.5000 & 100.0 & 1.0000 & 0.5000 & 100.0 & 1.0000 & 0.5000 \\
  &CLT[d4,2] & 98.0 & 0.4899 & 0.2685 & 92.2 & 0.5817 & 0.4966 & 94.4 & 0.4836 & 0.7228 \\
  &CI[d4,2] & 100.0 & 1.0000 & 0.5000 & 100.0 & 1.0000 & 0.5000 & 100.0 & 1.0000 & 0.5000 \\
  &CLT[d4,5] & 86.8 & 0.3477 & 0.2064 & 88.2 & 0.4848 & 0.4739 & 91.8 & 0.4564 & 0.7183 \\
  &CI[d4,5] & 100.0 & 1.0000 & 0.5000 & 100.0 & 1.0000 & 0.5000 & 100.0 & 1.0000 & 0.5000 \\ 
 \hline
  $n=100$&CLT[i2,2] & 97.0 & 0.4689 & 0.2628 & 94.0 & 0.5232 & 0.4939 & 98.0 & 0.4143 & 0.7604 \\
 & CI[i2,2] & 100.0 & 0.9997 & 0.4999 & 100.0 & 1.0000 & 0.5000 & 100.0 & 1.0000 & 0.5000 \\
 & CLT[i2,5] & 92.4 & 0.2907 & 0.1999 & 91.2 & 0.3670 & 0.4911 & 91.0 & 0.3521 & 0.7682 \\
 & CI[i2,5] & 100.0 & 0.9998 & 0.4999 & 100.0 & 0.9992 & 0.5004 & 100.0 & 0.9078 & 0.5461 \\
 & CLT[d4,2] & 97.6 & 0.3865 & 0.2299 & 93.6 & 0.4259 & 0.4900 & 93.8 & 0.3704 & 0.7690 \\
 & CI[d4,2] & 100.0 & 1.0000 & 0.5000 & 100.0 & 1.0000 & 0.5000 & 100.0 & 1.0000 & 0.5000 \\
 & CLT[d4,5] & 90.2 & 0.2691 & 0.1965 & 89.4 & 0.3509 & 0.4882 & 90.4 & 0.3486 & 0.7655 \\
 & CI[d4,5] & 100.0 & 1.0000 & 0.5000 & 100.0 & 0.9993 & 0.5003 & 100.0 & 0.9026 & 0.5487 \\
 \hline
$n=500$&  CLT[i2,2] & 95.8 & 0.2540 & 0.2057 & 92.8 & 0.2365 & 0.4997 & 94.0 & 0.2095 & 0.7983 \\
  &CI[i2,2] & 100.0 & 0.6990 & 0.3495 & 100.0 & 0.9399 & 0.5028 & 100.0 & 0.6864 & 0.6568 \\
  &CLT[i2,5] & 95.0 & 0.1363 & 0.2004 & 93.6 & 0.1657 & 0.4980 & 93.8 & 0.1712 & 0.7983 \\
  &CI[i2,5] & 100.0 & 0.5772 & 0.2886 & 100.0 & 0.7113 & 0.5192 & 100.0 & 0.5361 & 0.7319 \\
  &CLT[d4,2] & 95.2 & 0.1965 & 0.2032 & 93.8 & 0.1908 & 0.4987 & 94.2 & 0.1820 & 0.7982 \\
  &CI[d4,2] & 100.0 & 0.7002 & 0.3501 & 100.0 & 0.9459 & 0.5048 & 100.0 & 0.6806 & 0.6597 \\
  &CLT[d4,5] & 93.6 & 0.1250 & 0.1997 & 93.6 & 0.1586 & 0.4977 & 94.2 & 0.1700 & 0.7967 \\
  &CI[d4,5] & 100.0 & 0.5972 & 0.2986 & 100.0 & 0.7272 & 0.5316 & 100.0 & 0.5329 & 0.7335 \\
 \hline
 $n=1000$ &CLT[i2,2] & 95.4 & 0.1829 & 0.2019 & 93.8 & 0.1673 & 0.4988 & 94.4 & 0.1485 & 0.7988 \\
  &CI[i2,2] & 100.0 & 0.5500 & 0.2750 & 100.0 & 0.6912 & 0.5015 & 100.0 & 0.5441 & 0.7279 \\
  &CLT[i2,5] & 95.0 & 0.0963 & 0.1990 & 92.2 & 0.1173 & 0.4992 & 94.0 & 0.1211 & 0.7972 \\
  &CI[i2,5] & 100.0 & 0.4596 & 0.2302 & 100.0 & 0.5022 & 0.5092 & 100.0 & 0.4434 & 0.7779 \\
  &CLT[d4,2] & 94.6 & 0.1392 & 0.2009 & 93.2 & 0.1350 & 0.4981 & 93.8 & 0.1287 & 0.7979 \\
  &CI[d4,2] & 100.0 & 0.5491 & 0.2745 & 100.0 & 0.6873 & 0.5026 & 100.0 & 0.5412 & 0.7294 \\
  &CLT[d4,5] & 96.0 & 0.0884 & 0.1993 & 92.8 & 0.1123 & 0.4998 & 94.4 & 0.1203 & 0.7974 \\
  &CI[d4,5] & 100.0 & 0.4725 & 0.2365 & 100.0 & 0.5130 & 0.5168 & 100.0 & 0.4419 & 0.7790 \\
   \hline
$n=10000$&CLT[i2,2] & 95.0 & 0.0579 & 0.2001 & 95.2 & 0.0529 & 0.5010 & 95.4 & 0.0469 & 0.8007 \\
  &CI[i2,2] & 100.0 & 0.2179 & 0.2004 & 100.0 & 0.2179 & 0.5012 & 100.0 & 0.2179 & 0.8009 \\
  &CLT[i2,5] & 94.4 & 0.0305 & 0.2001 & 94.8 & 0.0371 & 0.5002 & 96.4 & 0.0383 & 0.8006 \\
  &CI[i2,5] & 100.0 & 0.1594 & 0.2008 & 100.0 & 0.1594 & 0.5009 & 100.0 & 0.1594 & 0.8013 \\
  &CLT[d4,2] & 95.0 & 0.0440 & 0.2001 & 95.2 & 0.0427 & 0.5006 & 95.6 & 0.0407 & 0.8007 \\
  &CI[d4,2] & 100.0 & 0.2165 & 0.2006 & 100.0 & 0.2165 & 0.5011 & 100.0 & 0.2165 & 0.8011 \\
  &CLT[d4,5] & 94.4 & 0.0280 & 0.2001 & 94.0 & 0.0355 & 0.5001 & 97.0 & 0.0381 & 0.8004 \\
  &CI[d4,5] & 100.0 & 0.1633 & 0.2020 & 100.0 & 0.1633 & 0.5020 & 100.0 & 0.1633 & 0.8023 \\
   \hline
\end{tabular}
\end{center}
\caption{Monte-carlo experiments based on $500$ replications of a fractional Brownian motion with Hurst parameter $H=0.2,0.5,0.8$ and scaling coefficient $C=1$ (assumed to be unknown), for $M=2,5$ and for different values of the sample size. The filters i2 and d4 denote respectively the filter of Increments of order 2 and the Daublets 4. For these simulations the vector $\Vect{d}$ has been fixed to the vector $\Vect{A}$.} \label{tab-simICgen}
\end{table}

\subsection{Discussion}
We propose non-asymptotic confidence intervals for the Hurst parameter of a standard or non-standard fBm based on concentration inequalities. 
They are computable in particular for small sample size and several theoretical improvements are obtained:
\begin{itemize}
\item When the scaling parameter $C$ is known, we have refined the confidence interval proposed in \cite{A-BreNouPec09}: 
the upper bound $H\leq H^\star<1$ is relaxed, the condition on the sample size $n$ is sharper 
and our new confidence intervals are valid for a large class of filter~$\ensuremath{{a^{}}}$.
\item As a by-product in our way to optimize the numeric bounds, we have slightly improved the bounds obtained by \cite{A-NouVie09} in the general concentration inequality (see Proposition~\ref{prop-CIH2}).
\item The case where $C$ is unknown has never been considered with concentration inequalities before Proposition \ref{prop-ICgen}.  
\item The asymptotic properties are similar to that of confidence intervals based on central limit theorems. 
More specifically, the length of the confidence intervals derived by concentration inequalities behaves asymptotically as the ones of confidence intervals based on central limit theorems, 
that is $1/(\sqrt{n}\log(n))$ when $C$ is known and $1/\sqrt{n}$ when $C$ is unknown.
\end{itemize}

The comparison with confidence interval based on CLT  is contrasted:  
while the Monte-Carlo experiments are correct when $C$ is known (in terms of coverage rate and of lengths of the confidence intervals), 
they are not good when $C$ is unknown: the lengths equal often $1$, {\it i.e.} the intervals correspond to $(0,1)$, when the sample size is small and are about five times larger when $n$ is large. 
In fact, the confidence intervals derived from concentration inequalities are too much "sympathetic": 
the coverage rate is rather far from $1-\alpha$ (based on $500$ replications, it is even often equal to $100\%$). 
From a statistical point of view, this is the main reason why the length of the confidence interval is sometimes much larger than the ones based on central limit theorems. \\
From a mathematical point of view, this is due to the fact that, in Proposition~\ref{prop-ICgen}, the dilatations of a filter are actually handled separately.
As a consequence, the errors induced by each dilatation, and controled by the concentration inequalities \eqref{eq-phir}--\eqref{def-phil}, add up, see \eqref{eq:addup}. 
This explains that the proposed confidence interval based on concentration inequalities are less performing in this case while, in comparison, multivariate CLT are used for standard confidence intervals. 
Improvements would require to use multivariate concentration inequalities, generalizing Proposition \ref{prop-CIH2}, which,  at the moment, are not available. 
This is the aim of future research to obtain such improvements.

As a conclusion, this work is the first attempt to define {\it computable} confidence intervals for the Hurst parameter $H$ 
of a standard and a non-standard fractional Brownian motion with another approach than the classical one based on central limit theorems 
(at the very exception of \cite{A-BreNouPec09} where the first non-asymptotic confidence intervals were derived for the standard fBM with a more theoretical motivation). 
We did not get around the question of the numerical performances via Monte-Carlo experiments. 
The conclusion is that, based on concentration inequalities, confidence intervals can be proposed for a large class of filters and without assumption on the Hurst parameter.  
The performances are comparable to the stantard confidence interval based on CLT when the scale parameter $C$ is known, 
while the procedure is underperforming when $C$ is unknown.  
This later case requires preliminary theoretical improvements for multivariate Gaussian quadratic forms that motivate our future studies.  



\appendix

\section{Exact computations of $\ell_1$-norm for filtered fBm}
\label{sec:exact}
In this section, we describe how explicit exact bound can be obtained for the correlation of a filtered fBm. 
Let $\ensuremath{{a}}$ be a filter of order $p$ and length $\ell$. 
Its covariance function is given by 
\begin{eqnarray*}
\pi^{\ensuremath{{a}}}_H(k)&=&-\frac 12 \sum_{q,r=0}^\ell a_qa_r|q-r+k|^{2H}
=-\frac 12\sum_{j=-\ell}^\ell \alpha_j|j+k|^{2H}
\end{eqnarray*}
where $\alpha_j=\sum_{\begin{subarray}{c}q,r=0\\ q-r=j\end{subarray}}^{\ell} a_qa_r$. 
Note that 
\begin{itemize}
\item $\alpha_j=\alpha_{-j}$, in particular $\pi^{\ensuremath{{a}}}_H(0)=-\sum_{j=1}^\ell\alpha_jj^{2H}$;
\item $\sum_{j=-\ell}^\ell \alpha_j=\sum_{q,r}^\ell a_qa_r=0$, 
\item for all $h\leq2p-1$, we have 
\begin{eqnarray}
\nonumber
\sum_{j=-\ell}^\ell j^h\alpha_j&=&\sum_{j=-\ell}^\ell j^h\sum_{q-r=j}a_qa_r
=\sum_{j=-\ell}^\ell \sum_{q-r=j}(q-r)^ha_qa_r
=\sum_{q,r=0}^\ell (q-r)^ha_qa_r\\
\nonumber
&=&\sum_{q,r=0}^\ell \sum_{k=0}^h {h\choose k}q^k(-r)^{h-k}a_qa_r \\
\nonumber&=&\sum_{k=0}^h \left((-1)^{h-k}{h\choose k}\left(\sum_{q=0}^\ell q^ka_q\right)\left( \sum_{q=0}^\ell r^{h-k}a_r\right)\right)\\
&=&0.\label{eq:alphaj}
\end{eqnarray}
\item $\sum_{j\not=0} \alpha_j=-\alpha_0=-\sum_{q=0}^\ell a_q^2<0$, $\alpha_\ell=a_0a_\ell$.
\end{itemize}
A crucial observation is that, at least for $|k|$ large enough, all the $\pi_H^\ensuremath{{a^{}}}(k)$, 
and thus all the $\rho^{\ensuremath{{a}}}(k)$, have the same sign. 
Indeed, using \eqref{eq:alphaj}, we have for $|k|\geq \ell$: 
\begin{eqnarray*}
\pi_H^\ensuremath{{a^{}}}(k)
&=& -\frac 12\sum_{j=1}^\ell \alpha_j\left(|k+j|^{2H}+|k-j|^{2H}-2|k|^{2H}\right)\\
&=& -\frac{|k|^{2H}}2\sum_{j=1}^\ell \alpha_j\left((1+j/k)^{2H}+(1-j/k)^{2H}-2\right)\\
&=& -|k|^{2H}\sum_{i=p}^{+\infty} \left(\frac{(2H)(2H-1)\dots(2H-2i+1)}{(2i)!k^{2i}}\left(\sum_{j=1}^\ell \alpha_jj^{2i}\right)\right)\\
&\sim& -|k|^{2H-2p} \frac{(2H)(2H-1)\dots(2H-2p+1)}{(2p)!}\left(\sum_{j=1}^\ell \alpha_jj^{2p}\right).
\end{eqnarray*}
This observation allows to reduce the computation of the $\ell^1$-norm $\|\rho_H^{\ensuremath{{a}}}\|_{\ell^1(\ZZ)}$, 
which is an infinite sum with modulus, to an infinite sum of correlations but without modulus plus some finite sum (with modulus remaining). 
Essentially, it remains to compute the sum of correlation without modulus. This is done below. 
But observe first that if there exists some $k(H,\ensuremath{{a}})\in\nit$ so that the correlations $\rho_H^{\ensuremath{{a}}}(k)$ have all the same sign for $|k|\geq k(H,\ensuremath{{a}})$ large enough. The value $k(H,\ensuremath{{a}})$ is not known in general. 
However for some family of filters (including increment-type filters $in$ and their dilatations $(in)^m$, $n,m\geq 1$), 
$k(H,\ensuremath{{a}})$ is known and explicit computations are tractable:
\begin{proposition}
\label{prop:conjecture-in}
For a dilated increment-type filter $\ensuremath{{a}}\in\{(in)^m\::\: n,m\geq 1\}$, we have
$k(H,\ensuremath{{a}})=\ell$, {\it i.e.} the following property holds true: 
\begin{equation}
\label{eq:conjecture}
\mbox{for all } |j|\geq \ell,\: \pi_H^{\ensuremath{{a}}}(j) \mbox{ is of the same sign as } (-1)^{p+1}(2H-1).
\end{equation}
\end{proposition}

\begin{demo}
Let $\theta_m(f)(x)=f(x+m)-2f(x)+f(x-m)$. 
Observe that if $f$ is a convex (resp. concave) function, then $\theta_m(f)(x)\geq 0$ (resp.  $\theta_m(f)(x)\leq 0$). 
For the $i1$ filter, we have $\pi_H^{i1}(x)=\frac 12\theta_1(|x|^{2H})$, 
for the $i2$ filter, we have $\pi_H^{i2}(x)=-\frac 12\theta_1^{\circ 2}(|x|^{2H})$ 
and more generally for the $m$-dilatation of the $in$ filter, we have $\pi_H^{(in)^m}(x)=\frac {(-1)^{n+1}}2\theta_m^{\circ n}(|x|^{2H})$.  

Observe also that the function $|x|^{2H}$ and all its iterated derivatives $(|x|^{2H})^{(2p)}$ of even order are convex if $H\geq 1/2$, concave if $H\leq 1/2$. 
By an immediate induction on $n$, we show that the same holds true for all $\theta_m^{\circ n}(|x|^{2H})$. 
In particular for $|j|\geq \ell m$, we obtain that $\pi_H^{(in)^m}(j)$ is of the same sign as $(-1)^{n+1}(2H-1)$. 
\CQFD
\end{demo}

\medskip
Obviously, the property \eqref{eq:conjecture} does not hold true for any filter (consider for instance $\{1,-4,5,-2\}$). In order to make easier our following explicit computation to derive exact value for  $\|\rho^{\ensuremath{{a}}}\|_{\ell^1(\ZZ)}$, 
we consider a filter $\ensuremath{{a}}$ satisfying \eqref{eq:conjecture} 
but we stress that for each particular filter the same strategy applies with some specific $k(H,\ensuremath{{a}})$.
First, for all $N\geq \ell$, we have: 
\begin{eqnarray*}
\nonumber
-2\sum_{j=\ell}^N \pi_H^{\ensuremath{{a}}}(j)&=&\sum_{j=\ell}^N \sum_{k=-\ell}^\ell \alpha_k|j+k|^{2H}\\
\nonumber
&=&\sum_{k=-\ell}^\ell \alpha_k\sum_{j=\ell}^N|j+k|^{2H}
=\sum_{k=-\ell}^\ell \alpha_k\sum_{j=\ell+k}^{N+k}|j|^{2H}\\
\nonumber
&=&\alpha_{-\ell}S_{N-\ell}^H+\sum_{k=-\ell+1}^\ell \alpha_k\left(S_{N+k}^H-S_{\ell+k-1}^H\right)\\
\nonumber
&=&\alpha_{-\ell}S_{N-\ell}^H+\sum_{k=-\ell+1}^\ell \alpha_k\left(S_{N-l}^H+\sum_{j=N-\ell+1}^{N+k}|j|^{2H}-S_{\ell+k-1}^H\right)\\
\nonumber
&=&\left(\sum_{k=-\ell}^\ell \alpha_k\right) S_{N-l}^H
+\left(\sum_{k=-\ell+1}^\ell \alpha_k \sum_{j=N-\ell+1}^{N+k}|j|^{2H}\right)
-\left(\sum_{k=-\ell+1}^\ell \alpha_k S_{\ell+k-1}^H\right)\\
\label{eq:xN}
&=&x_N-\sum_{k=-\ell+1}^\ell \alpha_k S_{\ell+k-1}^H
\end{eqnarray*}
where $S_k^H=\sum_{j=0}^k j^{2H}$ and
\begin{eqnarray}
\nonumber
x_N&=&\sum_{j=N-\ell+1}^{N+\ell}\left(|j|^{2H}\sum_{k=j-N}^\ell\alpha_k\right)\\
\nonumber
&=&|N+\ell|^{2H} \sum_{i=0}^{2\ell-1}\left(\left(1-\frac i{N+\ell}\right)^{2H}\sum_{k=\ell-i}^\ell\alpha_k\right)\\
\label{eq:DL}
&=&|N+\ell|^{2H} \sum_{i=0}^{2\ell-1}\left(\left(1-\frac {2Hi}{N+\ell}+\frac {2H(2H-1)i^2}{2(N+\ell)^2} 
+O\left(\frac 1{(N+\ell)^3}\right)\right)\sum_{k=\ell-i}^\ell\alpha_k\right)
\end{eqnarray}
But 
$$
\sum_{i=0}^{2\ell-1}\sum_{k=\ell-i}^\ell\alpha_k=\sum_{k=-\ell+1}^\ell(\ell+k)\alpha_k=\sum_{k=-\ell}^\ell(\ell+k)\alpha_k=0
$$
and 
\begin{eqnarray*}
\sum_{i=0}^{2\ell-1}\left(i\sum_{k=\ell-i}^\ell\alpha_k\right)
&=&\sum_{k=-\ell+1}^\ell\left(\alpha_k\sum_{i=\ell-k}^{2\ell-1}i\right)
=\sum_{k=-\ell}^\ell\left(\alpha_k\sum_{i=\ell-k}^{2\ell-1}i\right)\\
&=&\sum_{k=-\ell}^\ell\alpha_k\left(\frac{2\ell(2\ell-1)}{2}-\frac{(\ell-k)(\ell-k-1)}2\right)
=0
\end{eqnarray*}
because of \eqref{eq:alphaj}. 
We obtain $x_N=O\big((N+\ell)^{2H-2}\big) \to 0$, $N\to+\infty$.
Actually, expanding $(1-i/(N+\ell))^{2H}$ to the $(2p-1)$-th order in \eqref{eq:DL}, 
and since  $\sum_{i=1}^N i^k$ is a polynomial in $N$ of degree $k+1$, 
\eqref{eq:alphaj} shows that $x_N=O\big((N+\ell)^{2H-2p+1}\big)$.
Finally with the property \eqref{eq:conjecture}, we have: 
$$
2\sum_{j=\ell}^{+\infty} |\pi_H^{\ensuremath{{a}}}(j)|
=(-1)^{p+1}\epsilon((2H-1)\sum_{k=-\ell+1}^\ell \alpha_k S_{\ell+k-1}^H
$$
and 
\begin{eqnarray}
\nonumber
\|\rho^{\ensuremath{{a}}}_H\|_{\ell^1(\ZZ)}
&=&1+2\sum_{k=1}^{\ell-1}|\rho^{\ensuremath{{a}}}_H(k)|+2\sum_{k=\ell}^{+\infty}|\rho^{\ensuremath{{a}}}_H(k)|\\
\nonumber
&=&1
+\sum_{k=1}^{\ell-1}\left|\frac{\sum_{j=-\ell}^\ell \alpha_j|j+k|^{2H}}{\sum_{j=1}^{\ell}\alpha_j j^{2H}}\right|
+(-1)^{p+1}\epsilon(2H-1)\frac{\sum_{k=-\ell+1}^\ell \alpha_k S_{\ell+k-1}^H}{|\sum_{j=1}^\ell\alpha_j j^{2H}|}\\
\label{eq:normexplicite}
&=&1
+\sum_{k=1}^{\ell-1}\frac{\left|\sum_{j=-\ell}^\ell \alpha_j|j+k|^{2H}\right|}{-\sum_{j=1}^{\ell}\alpha_j j^{2H}}
+(-1)^{p+1}\epsilon(2H-1)\frac{\sum_{k=-\ell+1}^\ell \alpha_k S_{\ell+k-1}^H}{-\sum_{j=1}^\ell\alpha_j j^{2H}},
\end{eqnarray}
where we recall that $\epsilon(2H-1)=sign(2H-1)$. First, note that the modulus has been removed in the denominator of \eqref{eq:normexplicite} 
according to the following observation:
$$
\sum_{j=1}^\ell\alpha_j j^{2H}=\frac 12\sum_{j=-\ell}^\ell\alpha_j j^{2H}
\to_{H\to 0}\frac 12\sum_{j\not=0}\alpha_j=\sum_{j=-\ell}^\ell\alpha_j-\alpha_0=-\alpha_0<0.
$$
Since we assume moreover $\pi_H^{\ensuremath{{a}}}(0)\not=0$, this means that $\pi_H^{\ensuremath{{a}}}(0)>0$ 
and that $\left|\sum_{j=1}^\ell\alpha_j j^{2H}\right|=-\sum_{j=1}^\ell\alpha_j j^{2H}$.

\medskip
Next, note that \eqref{eq:normexplicite} is an explicit expression involving only finite sums and can be easily explicitely optimized for $H\in (0,1)$ for every given $\ensuremath{{a}}$ satisfying $\mathbf{H}^{\ensuremath{{a^{}}}}$.  
Note that, for $p\geq 2$, when $H\to 1$, right-hand side of \eqref{eq:normexplicite} remains well defined. 
Observe first that since for any fixed $k$, $\lim_{H\to 1}S_k^H=S_k^1=\frac{k(k+1)(2k-1)}{6}$, we have using \eqref{eq:alphaj}
$$
\lim_{H\to 1}\sum_{k=-\ell+1}^\ell \alpha_k S_{\ell+k-1}^H=
\frac 16\sum_{k=-\ell+1}^\ell \alpha_k (\ell+k-1)(\ell+k)(2\ell+2k-1)=0.
$$
The same holds true for $\sum_{j=-\ell}^\ell \alpha_j|j+k|^{2H}$ and $\sum_{j=1}^{\ell}\alpha_j j^{2H}$, 
but under $\mathbf{H}^{\ensuremath{{a^{}}}}$ in \eqref{eq-taua}, 
the rule of l'Hospital entails $\lim_{H\to 1^-} \|\rho^{\ensuremath{{a}}}_H\|_{\ell^1(\ZZ)}$ exists and is finite. 
Since obviously, $\|\rho^{\ensuremath{{a}}}_H\|_{\ell^1(\ZZ)}$ is a continuous function of $H\in [0,1)$, 
this ensures the continuity of 
$\|\rho^{\ensuremath{{a}}}_H\|_{\ell^1(\ZZ)}$ on $[0,1]$ 
and the constant $\kappa^{\ensuremath{{a}}}$ in our confidence interval is obtained by maximazing the explicit function in \eqref{eq:normexplicite}. 


\medskip
\noindent
{\bf Dilated simple increments $(i1)^m=\{-1,1\}^m$.} 
In this case, $\ell=m$, $p=1$, $\alpha_j=0$ for $1<j<m$ 
and $\alpha_0=2$, $\alpha_{\pm m}=-1$ so that \eqref{eq:normexplicite} rewrites:
\begin{equation}
\label{eq:i1H}
\left\|\rho_H^{\{-1,1\}^m}\right\|_{\ell^1(\ZZ)}=1+\sum_{j=1}^{m-1}\frac{\left||j+m|^{2H}-2|j|^{2H}+|j-m|^{2H}\right|}{m^{2H}}+\frac{S_{2m-1}^H-2S_{m-1}^H}{m^{2H}}.
\end{equation}
For instance for $m=1$, $\left\|\rho_H^{\{-1,1\}}\right\|_{\ell^1(\ZZ)}=2$ and for $m=2$, $\left\|\rho_H^{\{-1,1\}}\right\|_{\ell^1(\ZZ)}=2\frac{4^H+9^H-1}{4^H}$, 
so that $\kappa^{i1}=4$ and $\kappa^{(i1)^2}=8$ (recall that in this case, we optimize for $H\in(0,1/2]$). 

In general, since the right-hand side of \eqref{eq:i1H} is a continuous function of $H$, and since for all $k\geq 1$, $S_k^{1/2}=\frac{k(k+1)}2$, we have 
$\lim_{H\to (1/2)^-}\left\|\rho_H^{\{-1,1\}^m}\right\|_{\ell^1(\ZZ)}=2m$
while $\left\|\rho_{1/2}^{\{-1,1\}^m}\right\|_{\ell^1(\ZZ)}=m$, exhibiting a discontinuity of the $\ell^1$-norm for the dilated $i^1$ filters. 


\medskip
\noindent
{\bf Dilated double increments $(i2)^m=\{1,-2,1\}^m$.}
In this case, $\ell=2m$, $p=2$ and $\alpha_0=6$, $\alpha_{\pm m}=-4$, $\alpha_{\pm 2m}=1$, $\alpha_j=0$, $j\not=0,\pm m,\pm 2m$, 
so that \eqref{eq:normexplicite} rewrites:
\begin{eqnarray*}
\|\rho^{\{1,-2,1\}^m}_H\|_{\ell^1(\ZZ)}
&=&1+\sum_{k=1}^{2m-1}\frac{\left||k-2m|^{2H}-4|k-m|^{2H}+6|k|^{2H}-4|k+m|^{2H}+|k+2m|^{2H}\right|}{m^{2H}(4-4^H)}\\
&&+\epsilon(1-2H)\frac{-4S_{m-1}^H+6S_{2m-1}^H-4S_{3m-1}^H+S_{4m-1}^H}{m^{2H}(4-4^H)}.
\end{eqnarray*}
In order to obtain explicit values, we focus on the cases $m=1$ and $m=2$. 
First, for $m=1$, \eqref{eq:normexplicite}  reduces to 
\begin{eqnarray*}
\|\rho_H^{\ensuremath{{i2^{}}}}\|_{\ell^1(\ZZ)} 
&=&\left\{\begin{array}{ll}
1+ \frac{10-7\times 4^H + 2\times 9^H}{4-4^H},&H\leq 1/2\\
2,&H\geq 1/2
\end{array} \right.
\end{eqnarray*}
and elementary computations entail:
$$
\kappa^{i2}
= 2\times \lim_{H\to 0^+}\| \rho_H^{\ensuremath{{i2^{}}}}\|_{\ell^1(\ZZ)} 
= 2\left( 1+\frac{5}3 \right) =\frac{16}3.
$$
Next, for $m=2$, since  
$$\begin{array}{lllll}
2\pi_H^{(i2)^2}(1) &=& -2+3\times 9^H -25^H &\geq 0 &\forall H\in(0,1) \\
2\pi_H^{(i2)^2}(2) &=& -7\times 4^H + 4\times 16^H -36^H & \leq 0 &\forall H\in(0,1) \\
2\pi_H^{(i2)^2}(3) &=& 3-6\times 9^H + 4\times 25^H -49^H & \leq 0 & \forall H\in(0,1)
\end{array}
$$
expression \eqref{eq:normexplicite} reduces to 
\begin{eqnarray*}
\left\|\rho_H^{(i2)^2}\right\|_{\ell^1(\ZZ)} 
&=&\left\{ \begin{array}{ll}
1+ \frac{-6+10\times 4^H+12\times 9^H-7\times 16^H -8\times 25^H + 2\times 36^H + 2\times 49^H}{4^H(4-4^H)} & \mbox{ for } H\leq1/2 \\
1+ \frac{-4+4\times 4^H + 6\times 9^H -16^H -2\times 25^H}{4^H(4-4^H)} & \mbox{ for } H\geq1/2.
\end{array} \right.
\end{eqnarray*}
An elementary study of this function, together with the rule of l'Hospital, entails that 
$$
\kappa^{(i2)^2}
=2\times \sup_{H\in[0,1]} \left\| \rho_H^{(i2)^2}\right\|_{\ell^1(\ZZ)}  
= 2\times \lim_{H\to 1-}\left\| \rho_H^{(i2)^2}\right\|_{\ell^1(\ZZ)}
= 2\left(1 + \frac{25\log(5)-27\log(3)}{8\log(2)} \right)\simeq 7.813554.
$$






\bibliographystyle{plainnat.bst}

\bibliography{ic}

\begin{thebibliography}{13}
\providecommand{\natexlab}[1]{#1}
\providecommand{\url}[1]{\texttt{#1}}
\expandafter\ifx\csname urlstyle\endcsname\relax
  \providecommand{\doi}[1]{doi: #1}\else
  \providecommand{\doi}{doi: \begingroup \urlstyle{rm}\Url}\fi

\bibitem[Beran(1994)]{B-Ber94}
J.~Beran.
\newblock \emph{{Statistics for long-memory processes}}.
\newblock Chapman \& Hall/CRC, 1994.

\bibitem[Breton et~al.(2009)Breton, Nourdin, and Peccati]{A-BreNouPec09}
J.-C. Breton, I.~Nourdin, and G.~Peccati.
\newblock Exact confidence intervals for the hurst parameter of a fractional
  brownian motion.
\newblock \emph{Electron. J. Statist.}, 3:\penalty0 416--425, 2009.

\bibitem[Coeurjolly(2000)]{A-Coe00}
J.-F. Coeurjolly.
\newblock Simulation and identification of the fractional {Brownian} motion: a
  bibliographical and comparative study.
\newblock \emph{J. Stat. Softw.}, 5\penalty0 (7):\penalty0 1--53, November
  2000.

\bibitem[Coeurjolly(2001)]{A-Coe01}
J.-F. Coeurjolly.
\newblock Estimating the parameters of a fractional {Brownian} motion by
  discrete variations of its sample paths.
\newblock \emph{Stat. Infer. Stoch. Process.}, 4\penalty0 (2):\penalty0
  199--227, January 2001.

\bibitem[Coeurjolly(2008)]{A-Coe08}
J.-F. Coeurjolly.
\newblock Hurst exponent estimation of locally self-similar gaussian processes
  using sample quantiles.
\newblock \emph{J. Stat. Softw.}, 36\penalty0 (3):\penalty0 1404--1434, 2008.

\bibitem[Daubechies(2006)]{A-Dau06}
I.~Daubechies.
\newblock Orthonormal bases of compactly supported wavelets.
\newblock \emph{Communications on Pure and, Applied Mathematics}, 41\penalty0
  (7):\penalty0 909--996, 2006.

\bibitem[Doukhan et~al.(2003)Doukhan, Oppenheim, and Taqqu]{B-DouOppTaq03}
P.~Doukhan, G.~Oppenheim, and M.S. Taqqu.
\newblock \emph{{Theory and applications of long-range dependence}}.
\newblock Birkhauser, 2003.

\bibitem[Istas and Lang(1997)]{A-IstLan97}
J.~Istas and G.~Lang.
\newblock Quadratic variations and estimation of the h{\"o}lder index of a
  gaussian process.
\newblock \emph{Ann. Inst. H. Poincar{\'e} Probab. Statist.}, 33:\penalty0
  407--436, 1997.

\bibitem[Kent and Wood(1997)]{A-KenWoo97}
J.T. Kent and A.T.A. Wood.
\newblock Estimating the fractal dimension of a locally self-similar gaussian
  process using increments.
\newblock \emph{J. Roy. Statist. Soc. Ser. B}, 59:\penalty0 679--700, 1997.

\bibitem[Mandelbrot and Ness(1968)]{A-ManVan68}
B.~Mandelbrot and J.~Van Ness.
\newblock Fractional brownian motions, fractional noises and applications.
\newblock \emph{SIAM Rev.}, 10:\penalty0 422--437, 1968.

\bibitem[Nourdin and Viens(2009)]{A-NouVie09}
I.~Nourdin and F.G. Viens.
\newblock Density formula and concentration inequalities with malliavin
  calculus.
\newblock \emph{Elec. J. Probab.}, 14:\penalty0 2287--2309, 2009.

\bibitem[Percival and Walden(2000)]{B-PerWal00}
D.~B. Percival and A.~T. Walden.
\newblock \emph{Wavelet Methods for Time Series Analysis}.
\newblock Cambridge University Press, 2000.

\bibitem[Shen et~al.(2007)Shen, Zhu, and Lee]{A-SheZhuLee07}
Haipeng Shen, Zhengyuan Zhu, and Thomas C.~M. Lee.
\newblock Robust estimation of the self-similarity parameter in network traffic
  using wavelet transform.
\newblock \emph{Signal Process.}, 87\penalty0 (9):\penalty0 2111--2124, 2007.

\end{thebibliography}

\end{document}